\documentclass[aop,MSNbibl,seceqn,citesort,dvips]{arximspdf}

\doi{10.1214/12-AOP743} %kopijuoti is PTS
\volume{41}
\issue{3B}
\pubyear{2013}
\firstpage{1831}
\lastpage{1863}

\makeatletter
\newcommand{\eqref}[1]{(\ref{#1})}
\newcommand{\cF}{\mathcal{F}}
\newcommand{\cQ}{\mathcal{Q}}

\newcommand{\cD}{\mathcal{D}}
\newcommand{\cU}{\mathcal{U}}
\newcommand{\E}{\mathbb{E}}

\newcommand{\varqM}{\langle M \rangle}
\newcommand{\ind}{\mathbf{1}}

\newcommand{\varq}[1]{\langle #1 \rangle}
\newcommand{\expD}{\mathcal{D}_{\exp}}

\newcommand{\lo}[1]{\overline{#1}}

\newcommand{\lu}[1]{\underline{#1}}
\newcommand{\uq}{\underline{q}}
\newcommand{\urho}{\underline{\rho}}

\newtheorem{theorem}{Theorem}[section]
\newtheorem{elemma}[theorem]{Lemma}
\newtheorem{eproposition}[theorem]{Proposition}
\newtheorem{ecorollary}[theorem]{Corollary}

\newproclaim{edefinition}[theorem]{Definition}
\newproclaim{remark}{Remark}

\makeatother

\begin{document}
\begin{frontmatter}

\title{Monotone stability of quadratic semimartingales with
applications to unbounded general quadratic~BSDEs}
\runtitle{Monotone stability of quadratic semimartingales}

\begin{aug}
\author[A]{\fnms{Pauline} \snm{Barrieu}\corref{}\ead[label=e1]{p.m.barrieu@lse.ac.uk}}
\and
\author[B]{\fnms{Nicole} \snm{El Karoui}\ead[label=e2]{nicole.elkaroui@cmap.polytechnique.fr}\thanksref{t1}}
\thankstext{t1}{Supported in part by the ``Chaire Financial Risk'' of
the Risk Foundation, Paris.}
\runauthor{P. Barrieu and N. El Karoui}
\affiliation{London School of Economics and Universit\'{e} Pierre et
Marie Curie}
\address[A]{Statistics Department\\
London School of Economics\\
Houghton street\\
London WC2A2AE\\
United Kingdom\\
\printead{e1}} %adresu isvedimo komanda gale!
\address[B]{LPMA\\
Universit\'{e} Pierre et Marie Curie (Paris 6)\\
CNRS: UMR 7599\\
75005 Paris\\
France\\
\printead{e2}}
\end{aug}

% HISTORY:
\received{\smonth{1} \syear{2011}}
\revised{\smonth{1} \syear{2012}}

% ABSTRACT
%
\begin{abstract}
In this paper, we study the stability and convergence of some general
quadratic semimartingales. Motivated by financial applications, we
study simultaneously the semimartingale and its opposite. Their
characterization and integrability properties are obtained through some
useful exponential submartingale inequalities. Then, a general
stability result, including the strong convergence of the martingale
parts in various spaces ranging from $\mathbb{H}^1$ to BMO, is derived
under some mild integrability condition on the exponential of the
terminal value of the semimartingale. This can be applied in particular
to BSDE-like semimartingales.

This strong convergence result is then used to prove the existence of
solutions of general quadratic BSDEs under minimal exponential
integrability assumptions, relying on a regularization in both
linear-quadratic growth of the quadratic coefficient itself. On the
contrary to most of the existing literature, it does not involve the
seminal result of Kobylanski [\textit{Ann. Probab.} \textbf{28} (2010)
558--602] on bounded solutions.
\end{abstract}

% KEYWORDS
%


\end{frontmatter}

%====================================
%s1 #&#
\section{Introduction}\label{sec1}
%====================================
The Backward Stochastic Differential Equations (BSDEs) were first
introduced by Peng and Pardoux~\cite{PardouxPeng90} in 1990 in the
Lipschitz continuous framework, and then extended to continuous with
linear growth framework by Lepeltier and San Martin \cite
{LepeltierSanMartin97} in 1997. They have been soon recognized as powerful
tools with many different possible applications. More recently,
there has been an accrued interest for quadratic BSDEs, with various
fields of application such as risk sensitive control problems or
dynamic financial risk measures and indifference pricing in
mathematical finance.

In this case, the BSDE is an
equation of the following type:
%
%e1 #&#
%
\begin{equation}
-dY_{t}=g(t,Y_{t},Z_{t})\,dt-Z_{t}\,dW_{t},\qquad Y_{T}=\xi_{T},
\end{equation}
where $W_{\cdot}$ is a standard Brownian motion, and the
coefficient $g$ satisfies the following quadratic \textit{structure
condition} ${\mathcal Q}(l_{\cdot},c_{\cdot},\delta)$:
%
%e2 #&#
%
\begin{equation}%\label{eqquadraticgrowth}
| g(t,y,z)|\leq\kappa(t,y,z) \equiv\frac{1}{\delta
}l_{t}+c_t|y|+\frac
{\delta}{2}
|z|^{2}, \qquad d\mathbb{P}\otimes dt\mathrm{\mbox{-}a.s.},
\end{equation}
where $\delta>0$ is a given constant, and $(l_t), (c_t)$ are
predictable nonnegative processes.

The first result concerning the existence and uniqueness of solutions
to these equations was obtained in the bounded case in a Brownian
filtration setting by Kobylanski~\cite{Kobylanski00} in 2000. The proof
first relies on an exponential transformation as to come back to the better
known framework of BSDEs with a coefficient with linear growth and
then uses a regularization procedure to take the limit. The major
difficulty is then about proving the strong convergence of the
martingale parts without having to impose too strong
assumptions. This seminal paper has been extended in several
directions, to a continuous setting by Morlais
\cite{Morlais1}, to unbounded solutions by Briand
and Hu~\cite{Briand-Hu} or more recently by Mochel and Westray~\cite
{Mocha-Westray}. Some other authors have obtained
further results in some particular situations (see, e.g., Hu
and Schweizer~\cite{Hu-Schweizer},
Hu, Imkeller and Muller~\cite{Hu-Imkeller-Muller}, Mania and Tevzadze
\cite{ManiaTevzadze} or Delbaen, Hu and Richou \cite
{Delbaen-Hu-Richou}). Recently in 2008, Tevzadze~\cite{Tevzadze} has
given a direct proof for the existence and uniqueness of a bounded
solution in the Lipschitz-quadratic case.

%A lot of extensions are provi
%The question of existence and uniqueness of solutions to these
%equations with quadratic growth was first examined by Kobylanski
%setting, and
%was then extended to a continuous filtration setting by Morlais
%given a direct proof for the existence and uniqueness of a bounded
%solution in the Lipschitz-quadratic case. Briand
%and Hu (\cite{Briand-Hu} and~\cite{Briand-Huconvex}) have extended the
%existence result to
%unbounded solution and proved uniqueness for a convex coefficient.
%Some other authors have obtained
%further results in some particular situations (e.g. Hu and
%%Schweizer~\cite{Hu-Schweizer}, Hu, Imkeller and Muller
%In general, the proof relies
%first on an exponential transformation as to come back to the better
%known framework of BSDEs with a coefficient with linear growth and
%then uses a regularization procedure to take the limit. The major
%difficulty is then about proving the strong convergence of the
%martingale parts without having to impose too strong
%assumptions.\vspace{2mm}\\
%Keeping in mind a possible application {\cgblue to BSDEs},
We adopt in this paper a
completely different approach and consider a forward point of view to
treat directly the questions of convergence.
% of {\cgblue the stability of a such class of semimartingales, in
%particular the strong convergence of the martingale parts.}
%
To do so, we introduce the notion of general quadratic semimartingales
in Section~\ref{sectionquadraticsemimartingale} and study their
characterization with regards to their
integrability properties under some interesting exponential
transformations in Section~\ref{subsectioncharactexpoinequality}.
Mainly motivated by
financial applications, where a seller price and a buyer price have to
be given simultaneously, we apply systematically the same
assumptions on the semimartingale and on its opposite. Having both
exponential integrability properties proves to be essential in the a
priori estimation of their quadratic variations. In Section \ref
{sectionestimates}, we
obtain a general stability result, including the strong convergence of
the martingale parts as presented in Theorem~\ref{thNicolas}. The
result is very general and simply require the existence of exponential
moment of the absolute value of (or quantities related to) the terminal
value of the semimartingales. Our approach allows us to obtain the
strong convergence of the martingale parts in $\mathbb{H}^1$. Stability
results are also obtained in various spaces, depending on the
assumption made on the terminal values. It is interesting to note that,
on the contrary to most of the existing literature, the space of BMO
martingales does not play any particular role\vadjust{\goodbreak} as the semimartingales
are no longer bounded. This stability result is completed, in the BSDE
framework, by the convergence in total variation of the finite
variation part. In Section~\ref{sec5}, existence results become a possible
application of this stability result.
%BSDEs
More precisely, coming back to our initial motivation of
quadratic BSDEs, we first regularize the quadratic coefficient of
the BSDE through inf-convolution as to transform it into a
coefficient with linear-quadratic growth. This regularization as
linear-quadratic, and not simply linear, allows us to consider
situations which are typically not considered in the literature.
Applying the stability result of the previous section, we can pass to
the limit and prove the existence result for general quadratic BSDEs,
under ``minimal'' integrability assumptions. The power of the forward
point of view is striking as existence results are easily obtained in a
more general framework than the classical existing literature. However,
uniqueness results requires stronger assumptions on the solutions, as
in Kobylanski~\cite{Kobylanski00} for the bounded case, or for convex
BSDEs, as in Briand and Hu~\cite{Briand-Hu} or more recently in Mochel
and Westray~\cite{Mocha-Westray} with exponential moments of any order,
or in Delbaen, Hu and Richou~\cite{Delbaen-Hu-Richou} under weaker
integrability assumptions.

This approach has also other potential applications that we will not
discuss here for lack of space. We can just mention numerical
simulations of quadratic BSDEs, study in terms of risk measures and
dual representation, solving of associated HJB-type equations.
%The paper is organized as follows. In Sections 2 and 3, we introduce
%quadratic semimartingales, satisfying a quadratic structure
%condition similar to (\ref{eqquadraticgrowth}), derive some key
%properties, discuss their integrability and obtain their
%characterization using exponential inequalities. Section 4 is
%dedicated to the question of stability and convergence of the
%quadratic semimartingales. Quadratic variation estimates are first
%obtained, and then a general stability result for the quadratic
%semimartingales is presented in Theorem~\ref{thNicolas}. {\cgreen NEK:
%la fin est a preciser mieux} \textit{Finally,
%Section 5 comes back to the study of general quadratic BSDEs and in
%particular the existence of a solution in the light of the previous
%forward results of this paper.}
%====================================
%s2 #&#
\section{Quadratic semimartingales}\label{sectionquadraticsemimartingale}
%====================================
Quadratic BSDEs have recently received a lot of attention, mainly
due to the wide range of possible applications, involving
optimization problems with an exponential criterion, such as
risk-sensitive control problems introduced by Fleming in the 1980s (see
Fleming and Sheu~\cite{FlemingSheu} for financial applications, or El
Karoui and Hamad\`{e}ne for an application to risk-sensitive zero-sum
stochastic functional games~\cite{ElKaroui-Hamadene}).

Financial applications have generated a renewed interest for this type
of BSDEs, particularly in connection with the theory of dynamic risk
measures as in Barrieu and El Karoui~\cite{Barrieu-ElKaroui6}, or
indifference pricing with exponential utility (see, e.g., Rouge
and El Karoui~\cite{Rouge-ElKaroui}, Mania and Schweizer \cite
{Mania-Schweizer} or the recent book edited by Carmona~\cite{Carmona}
among many other
references). Therefore, it is particularly relevant to understand the
structure of these processes, and to obtain conditions ensuring their
stability.

In the classical martingale theory, Burkolder--Davis--Gundy-type
estimates are crucial to obtain convergence results for martingales in
$\mathbb{H}^p$ from the convergence of their terminal values. The study
of classical BSDEs with linear growth relies also on precise a priori
estimates coming from the martingale theory, arising from a forward
point of view (see, e.g., in a general framework, El Karoui and
Huang~\cite{NEKHuang}).
In this section, after having defined quadratic BSDEs,
we adopt a forward point of view, introducing quadratic
semimartingales, with a similar structure condition, studying their
main properties and deriving some characterization results, which
depend on various integrability assumptions. These results will be
very useful to derive some stability and convergence results in the
next section.
%====================================
%s2.1 #&#
\subsection{Definition of quadratic BSDEs and quadratic semimartingales}\label{paragraphdefinitionofquadraticsemimartingales}
%====================================
Let us briefly recall the definition of a quadratic BSDE. Let $(\Omega,
\mathcal{F},\mathbb{ P},(\mathcal{F}_{t}))$ be a filtered probability
space, where the filtration $(\mathcal{F}_{t})$
satisfies the usual conditions of completeness and right-continuity.
The $\sigma$-field on $\Omega\times\mathbb R^+$ generated by the
adapted and left continuous processes is called the predictable $\sigma
$-field and denoted by $\mathcal{P}$. In this paper, we only consider
\textit{continuous} filtered probability space, that is, a filtered
probability space such that any locally bounded martingale is a
continuous martingale. A classical example is the probability space
generated by a Brownian motion, and satisfying the usual conditions.
%
%================================
%pa2.1.0.1 #&#
\subsubsection*{Definition of quadratic BSDEs}
A~quadratic BSDE is an equation of the following type:
%
%e3 #&#
%
\begin{equation}
-dY_{t}=g(t,Y_{t},Z_{t})\,dt-Z_{t}\,dW_{t},\qquad Y_{T}=\xi_{T},
\label{eqBSDE}
\end{equation}
where $T>0$ is a given time horizon (possibly $\mathcal
{F}_{t}$-stopping time), $W_{\cdot}$ is a standard $d$-dimensional
$(\mathbb
{P},(\mathcal{F}_{t}))$-Brownian motion, and $Z_{t} \,dW_{t}$
simply denotes the scalar product. The $\mathcal{F}_{T}$-random
variable $\xi_{T}$ is the terminal condition,\setcounter{footnote}{1}\footnote{As pointed out
by one referee, the random variable $\xi_T$ has to be in fact
$\mathcal
{F}_{T^{-}}$-measurable as terminal value of a continuous process.}
and the coefficient $g$ is a $\mathcal{P}\otimes\mathcal{B}(\mathbb
{R}\times\mathbb{R}^d)$ measurable process satisfying the following
quadratic \textit{structure condition} ${\mathcal
Q}(l,c,\delta)$:
%
%e4 #&#
%
\begin{equation}\label{eqquadraticgrowth}
| g(\cdot,t,y,z)|\leq\kappa(t,y,z) \equiv|l_{t}|+c_t|y|+\frac
{\delta}{2}
|z|^{2}, \qquad d\mathbb{P}\otimes dt\mathrm{\mbox{-}a.s.},
\end{equation}
where $\delta>0$ is a given constant, and $(l_{\cdot}), (c_{\cdot})$ are
predictable positive\footnote{In the rest of the paper, we adopt the
following European terminology: a positive random variable $X$ verifies
$\mathbb P(X\geq0)=1$, and a strictly positive random variable
verifies $\mathbb P(X> 0)=1$. In the same way, a c\`{a}dl\`{a}g process
$K_{\cdot}$ is said to be increasing when, for any~$t,s$ such that $t
\geq s$, the random variable $A_t-A_s$ is positive and strictly increasing
when $A_t -A_s$ is strictly positive.} processes.

By solution to the $\operatorname{BSDE}(g,\xi_{T})$ defined in equation
(\ref
{eqBSDE}), we mean a pair of predictable processes taking values in
$\mathbb{R}\times\mathbb{R}^d$, $(Y,Z)=\{(Y_t,Z_t); t\in[0,T]\}$,
such that the paths of $Y$ are continuous,
$\int_0^T |Z_t|^2\,dt<\infty$, $\int_0^T
|g(t,Y_t,Z_t)|\,dt<\infty$, $\mathbb{P}$-a.s., and
%
%e5 #&#
%
\begin{equation}
\label{eqdefbsde}
Y_t=\xi_{T}+\int_t^T
g(s,Y_s,Z_s)\,ds-\int_t^T Z_s \,dW_s, \qquad\mathbb{P}\mathrm{\mbox
{-}a.s.}
\end{equation}
Note that, in the rest of the paper, this type of equality between two
processes has to be understood as holding up to
indistinguishability.\vadjust{\goodbreak}

This minimal definition will be completed later on by some further
integrability assumptions.
% allowing us to obtain some stability results and some conditions for
%the existence of a solution.
%
%pa2.1.0.2 #&#
\subsubsection*{Definition of quadratic semimartingales}
Adopting a forward point of view, a solution of a quadratic BSDE is a
quadratic It\^o's semimartingale~$Y_{\cdot}$, where the predictable
process with finite variation satisfies the same quadratic structure
condition \eqref{eqquadraticgrowth}. Such a condition needs to be
further specified when considering the more general framework of
quadratic semimartingales defined on a continuous filtered probability
space.
%
%de2.1 #&#
%
\begin{edefinition}[(Quadratic semimartingale)]\label
{defquadraticsemimartingale}
Let $Y_{\cdot}$ be a continuous semimartingale,
%defined on a continuous filtered probability space,
with the decomposition
$Y_{\cdot}=Y_0-V_{\cdot}+M_{\cdot}$, where $V_{\cdot}$ is a
predictable process with finite total
variation~$|V|_{\cdot}$ and $M_{\cdot}$ is a local martingale with
quadratic variation
$\varqM_{\cdot}$.

$Y_{\cdot}$ is a \textit{quadratic semimartingale} if there exist two adapted
continuous increasing
processes $\Lambda_\cdot$ and $C_\cdot$ and a positive constant
$\delta$, such
that the
structure condition $\mathcal{Q}(\Lambda,C,\delta)$ holds true:
%e6 #&#
%
\begin{equation}\label{eqstructurecondition}
d|V|_t \ll\frac{1}{\delta}\,d\Lambda_t+|Y_t| \,dC_t+\frac{\delta
}{2}\,d\varqM_t,\qquad
d\mathbb{P}\mathrm{\mbox{-}a.s.}
\end{equation}
The symbol $\ll$ stands for the strong order of increasing
processes, stating that the difference is an increasing process.
Sometimes we use the short notation $D^{\Lambda,C}_{\cdot}(Y,\delta
)=\frac
{1}{\delta}\Lambda_{\cdot}+|Y_{\cdot}|\star C_{\cdot}$, and even
simply $D^{\Lambda,C}_{\cdot}$
when there is no ambiguity.
At this stage, no particular integrability assumption is made on the
processes~$\Lambda_{\cdot}$ and $C_{\cdot}$.
\end{edefinition}

%At this stage, no particular integrability assumption is made on the
%processes $\Lambda_{\cdot}$ and $C_{\cdot}$.\\

\textit{Comments}:
(i) Observe that if $Y_{\cdot}$ is a quadratic semimartingale,
then $-Y_{\cdot}$ is also a quadratic semimartingale.\vspace*{-6pt}

\begin{longlist}[(iii)]
\item[(ii)]More generally, if $Y_{\cdot}$ is a quadratic semimartingale and
$\delta>0$, $Y^\delta_{\cdot} \equiv\delta Y_{\cdot}$ is a semimartingale
associated with $M^\delta_{\cdot} \equiv\delta M_{\cdot}$ with
quadratic variation
$\varq{M^\delta}_{\cdot}=\delta^2 \varq{M}_{\cdot}$ and $V^\delta
_{\cdot} \equiv\delta
V$. Then the structure condition for the process $\delta Y_{\cdot}$ becomes
$d|V^\delta|_t \ll \,d\Lambda_t+|Y^\delta_t| \,dC_t+\frac{1}{2}\,
d\varq
{M^\delta}_t$. This property justifies our choice of restricting our
study to quadratic semimartingales with constant $\delta=1$, without
any loss of generality.

\item[(iii)]The following notation specify different classes of quadratic
semimartingales, $\mathcal{Q}(\Lambda,C,\delta)$ for the general case,
$\mathcal{Q}(\Lambda,C)$ when $\delta=1$, $\mathcal{Q}$ when
$\Lambda
_{\cdot}\equiv0, C_{\cdot}\equiv0, \delta=1.$
\end{longlist}
%
%s2.2 #&#
\subsection{Exponential transformations and algebraic characterization of quadratic semimartingales}\label{parrecallsemimartingale}\label{parbasicproperties}%\label{parbasicproperties}
%%%%=============================================================================
%pa2.2.0.1 #&#
\subsubsection*{Some recalls on semimartingales on a continuous
probability space} (i) Let us first recall the conventional notation
for the exponential martingale of a continuous (local) martingale
$M_{\cdot}$ with quadratic variation $\varqM_{\cdot}$
%e7 #&#
%
\begin{equation}\label{eqexpmartingale}
\mathcal E_{\cdot}(M) \equiv\exp\bigl(M_{\cdot}- \tfrac{1}{2}\varqM
_{\cdot}\bigr).\vadjust{\goodbreak}
\end{equation}

\begin{longlist}[(iii)]
\item[(ii)]
% A \begin{edefinition}\label{defgeneralsubmartingale}
A right continuous left limited submartingale (c\`{a}dl\`{a}g in the
French denomination) $S_{\cdot}$ is a c\`{a}dl\`{a}g optional process
$S_{\cdot}=S_0+N_{\cdot}+K_{\cdot},$ where $N_{\cdot}$ is a local
martingale and $K_{\cdot}$ a
predictable c\`{a}dl\`{a}g increasing process.
The pair $(N_{\cdot},K_{\cdot})$ is called the additive decomposition
of $S$. When
$S_{\cdot}$ is a positive submartingale, $(M_{\cdot},A_{\cdot})$ is
said to
be the multiplicative decomposition of $S_{\cdot}$ if $S_{\cdot}=S_0
\mathcal{E}_{\cdot}(M)\exp(A_{\cdot})$, where~$M_{\cdot}$ is a
local martingale and
$A_{\cdot}$ a predictable c\`{a}dl\`{a}g increasing process.

\item[(iii)] Dellacherie and Meyer~\cite{Dellacherie-Meyer} (in Appendix
1---Probabilit\'{e}s et Potentiel~B) have extended this definition to right
and left limited submartingales (also known as strong submartingales)
when the increasing predictable process $K_{\cdot}$ is only with left and
right limits (l\`{a}dl\`{a}g in the French denomination), with the
following decomposition $K_{\cdot}=K^1_{\cdot}+K^2_{-\cdot}$, where
$K^1_{\cdot}$ is a c\`
{a}dl\`{a}g predictable increasing process and $K^2_{-\cdot}$ is the
process of the left limits of a c\`{a}dl\`{a}g optional increasing
process $K^2_{\cdot}$.
\end{longlist}

%
%==================================
%pa2.2.0.2 #&#
\subsubsection*{Characterization of $\mathcal{Q}$-semimartingales when
$\Lambda\equiv C \equiv0$ and $\delta=1$}
%
%For the sake of clarity of the presentation, we make a distinction
%between the $\mathcal{Q}$-semimartingales for which the
%processes $C_{\cdot}$ and $\Lambda_{\cdot}$ are identically equal to
%$0$ and $
%the general $\mathcal{Q}(\Lambda,C,\delta)$-semimartingales. This
%distinction appears naturally when considering dynamic risk measures
%(see e.g. %Barrieu and El Karoui~\cite{Barrieu-ElKaroui6}) %or
%Delbaen, Hu and Bao~\cite{Delbaen-Hu-Bao})
%and allows to derive very simple results.\\
The simplest $\mathcal{Q}$-semimartingales are those for which the
structure condition ${\mathcal Q}$ is saturated, that is, $V_{\cdot
}=\frac
{1}{2}\varqM_{\cdot}$ or $\lu{V}_{\cdot}=-\frac{1}{2}
\varqM_{\cdot}$. Because of their importance, we refer to them as $q
$ (resp.,
$\lu{q}$) semimartingales, and denote them by
%Because of their importance, we give them a specific
%denomination and refer to them as $q $ (resp. $\lu{q}$)
%semimartingales.
%$q_\delta$-(resp.
%$\lu{q}_\delta$-) semimartingales. In particular, when $\delta=1$
%(resp. $\delta=-1$), $q $ (resp. $\lu{q}$) semimartingales are also
%
%e8 #&#
%
\begin{equation}\label{eqrtM}
\cases{
r_{\cdot}(r_0,M) \equiv r_0+M_{\cdot}-\frac{1}{2}\varqM_{\cdot}
\equiv
r_0+r_{\cdot}(M),\vspace*{2pt}\cr
\lu{r}_{\cdot}(r_0,M) \equiv\lu{r}_0+ M_{\cdot}+\frac
{1}{2}\varqM
_{\cdot} \equiv\lu{r}_0-r_{\cdot}(-M).}
\end{equation}
The operator $M \rightarrow r_{\cdot}(M)$ is not an additive operator,
nevertheless $r_{\cdot}(M)+r_{\cdot}(M')=r_{\cdot}(M+M')+\langle
M,M' \rangle_{\cdot}$ and
$r_{\cdot}(M)-r_{\cdot}(M')=r_{\cdot}(M-M')-\langle M-M', M' \rangle
_{\cdot}$.

Taking the exponential of $r_{\cdot}(M)$ immediately leads to the
exponential martingale $\mathcal E_{\cdot}(M)=e^{r_{\cdot}(M)}$
defined in \eqref
{eqexpmartingale}, whilst the exponential of $\lu{r}_{\cdot}(M)$
leads to
$e^{\lu{r}_{\cdot}(M)}=(\mathcal E_{\cdot}(- M))^{-1}$.

%E(M)_{\cdot}\\
% \end{array}
% .
%Observe that $e^{r_{\cdot}(M)}$ is the exponential martingale $\mathcal
%E(M)=\exp(M- \demi\varqM_{\cdot} $, and
%$e^{\lu{r}_{\cdot}(\lu M)}=(\mathcal E(-\lu M))^{-1}$.}
%
It will also be interesting to introduce some asymmetry in the previous
definition of $\mathcal{Q}$-semimartingales, with the notion of
$\mathcal{Q}$-submartingales, especially useful when characterizing
the former.
%de2.2 #&#
%
\begin{edefinition}\label{defquadraticsubmartingale}
A $\cQ$-submartingale is a continuous (or l\`{a}dl\`{a}g) semimartingale
$X_{\cdot}=X_0-V_{\cdot}+M_{\cdot}$ such that $A_{\cdot}  \equiv-V_{\cdot
}+\frac{1}{2}\varqM_{\cdot}$ is a
predictable increasing process.
Equivalently, $e^{X_{\cdot}}=e^{X_0+A_{\cdot}}\mathcal E_{\cdot}(M)$
is a continuous (l\`{a}dl\`{a}g) submartingale.
%$X_{\cdot}=X_0+r_{\cdot}(M)+A_{\cdot}$ and
%only if $Y_{\cdot}$ and $-Y_{\cdot}$ are $\cQ$-submartingales.}
\end{edefinition}

Obviously a $\cQ$-semimartingale is a $\cQ$-submartingale. Remarkably,
applying this property to both $X$ and $-X$ is sufficient to
characterize $\cQ$-semimartingales. From a financial point of view,
this means that the same rules have to be used to characterize both the
buyer's and the seller's price.

%th2.3 #&#
%
\begin{theorem}\label{thcharacterizationcQquasimart}
Let $X_{\cdot}$ be a l\`{a}dl\`{a}g optional process. Then, $X_{\cdot
}$ is a $\cQ
$-semimartingale if and only if both processes $X$ and $-X$ are $\cQ
$-submartingales, or equivalently if and only if $\exp(X_{\cdot})$
and $\exp(-X_{\cdot})$
are submartingales. In all cases, $X_{\cdot}$ is a continuous process.
\end{theorem}

\begin{pf} We only have to prove the sufficiency.
% if $Y_{\cdot}$ is a $\cQ$-semimartingale, both processes $Y_{\cdot}$
%and $-Y_{\cdot}$
%are $\cQ$-submartingales, and then their exponential %$\exp(Y_{
Assume that $\exp(X_{\cdot})$ and $\exp(-X_{\cdot})$ are two
l\`{a}dl\`{a}g submartingales, with respective multiplicative
decomposition $(\lo{ M}_{\cdot},\lo A_{\cdot})$, and
$(\lu{ M}_{\cdot},\lu A_{\cdot})$. Taking the logarithm leads to two
different
decompositions of $X$,
\[
%which are compatible only if X is continuous.
%Indeed, from the multiplicative submartingale decomposition, we have:
X_{\cdot}= X_0+\lo{ M}_{\cdot}-\tfrac{1}{2}\varq{\lo
M}_{\cdot}
+ \lo A_{\cdot}\quad \mbox{and} \quad {-}X_{\cdot}= -X_0+
\lu{M}_{\cdot}-\tfrac{1}{2}\varq{\lu{M}}_{\cdot} +\lu{A}_{\cdot}.
\]
Since the martingales and their quadratic variations are continuous,
the jumps of~$X$ are the same as the positive jumps
of the increasing process $\lo A_{\cdot}$. The same remark holds true
for the
jumps of the process $-X$. As, the jumps of $X$ are simultaneously
positive and negative, the process $X_{\cdot}$ is continuous.

Moreover, from the uniqueness of the predictable decomposition of
$X_{\cdot}$ we know that $\lu{M}_{\cdot}=-\lo M_{\cdot}$. Hence,
$\varq{\lu{M}}=\varq{\lo{M}}$
and $\lo{A}_{\cdot}+\lu{A}_{\cdot}=\varq{M}_{\cdot}$.
From Radon--Nikodym's theorem, there exists a predictable process
$\alpha_{\cdot}$, with $0\leq\alpha_{t}\leq2$, such that
$d\lo A_{t}=\frac{1}{2}\alpha_{t}\,d \varq{M}_t$. Substituting $\lo
A_{\cdot
} $
into the decomposition of $X_{\cdot}$, we get $dX_{t}=-\frac
{1}{2}(1-\alpha
_{t})\,d\varq{M}_t+dM_t$
with $|1-\alpha_{t} |\leq1$. Therefore, $X_{\cdot}$ is a ${\mathcal
Q}$-semimartingale.
\end{pf}
%
%================================================
%pa2.2.0.3 #&#
\subsubsection*{Characterization of $\mathcal{Q}(\Lambda,C,\delta
)$-semimartingales via exponential transformation}
%Note that, if $Y_{\cdot}$ is a $\mathcal{Q}$-semimartingale, then both
%$Y_{\cdot}$
%and $-Y_{\cdot}$ are $\cQ$-submartingales.\vspace{1mm}\\
%After this brief introduction of quadratic semimartingales and some
%related notions, let us now
%look at some key properties.
In the general structure condition
(\ref{eqstructurecondition}), the presence of the term
$|Y_{\cdot}|\star C_{\cdot}$ makes the characterization of quadratic
semimartingales more difficult to obtain.
% than for $\mathcal{Q}$-semimartingales.
Nevertheless the transformations
proposed in the following proposition can partially reduce the problem to
$\cQ$-submartingales.
% \\[-12mm]
%
%th2.4 #&#
%
\begin{theorem}\label{thcharacterizationgeneralquasimart}
Let us introduce the following transformations of any adapted (l\`{a}dl\`{a}g) process $Y_{\cdot}$:
%e9 #&#
%e10 #&#
%
\begin{eqnarray}\label{eqXLambda,C}
X^{\Lambda,C}_t(Y) &\equiv& Y_t +\Lambda_t+ \int_0^t |Y_s|\,dC_s
\equiv
Y_t +D^{\Lambda,C}_t(Y), \label{eqXLambda,C}\\
 U^{\Lambda,C}_t(e^Y)& \equiv& e^{Y_t}+\int_0^te^{Y_s}\,d\Lambda
_s+\int
_0^te^{Y_s}|Y_s|\,dC_s. \label{eqULambda,C}
\end{eqnarray}
Then, $Y_{\cdot}$ is a $\mathcal{Q}(\Lambda,C,\delta
)$-semimartingale if and
only if $X^{\Lambda,C}_{\cdot}(\delta Y)$ and\break $X^{\Lambda,C}_{\cdot
}(-\delta Y)$
are $\cQ$-submartingales, or equivalently if and only if both processes
$U^{\Lambda,C}_{\cdot}(e^{\delta Y})$ and $U^{\Lambda,C}_{\cdot
}(e^{-\delta Y})$
are submartingales.
%la meme caracterisation est vraie en utilisant $U^{\Lambda,C}_t(e^Y)$
%a la place de $X^{\Lambda,C}_{\cdot}(Y)$.
%are submartingales.
\end{theorem}

The link between the two transformations $X^{\Lambda,C}$ and\vspace*{1pt}
$U^{\Lambda
,C}$ is clear when $Y$ is a continuous semimartingale, since
$dU^{\Lambda,C}_t(e^{Y})=e^{-D^{\Lambda,C}_t}\,de^{X^{\Lambda,C}_t(Y)}$
(see proof below). The motivation behind the transformation $U^{\Lambda
,C}_t(e^Y)$, first introduced by Briand and Hu~\cite{Briand-Hu} will be
presented later in Section~\ref{subsectioncharactexpoinequality}.\vadjust{\goodbreak}
%when class $(\mathcal{D})$ integrability conditions will be
%considered. In this case, such a transformation appears to be
%essential to obtain useful inequalities for %the derivation of some
%stability and convergence results.
%
%As we will see in details, this exponential-type transformation will
%prove t%
%
\begin{pf*}{Proof of Theorem~\ref{thcharacterizationgeneralquasimart}} We can assume $\delta=1$ without any loss of generality
[refer in particular to Comment (ii) at the end of Section \ref
{paragraphdefinitionofquadraticsemimartingales}].
\begin{longlist}[(ii.b)]
\item[(i.a)]Necessary condition:
%Since $\delta Y$ is a $\mathcal{Q}(\Lambda, C)$ semimartingale, it is
%sufficient to study the case $\delta=1$.
Let $\alpha^V_{\cdot}\in[-1,1]$ be a predictable process such that
$V_{\cdot}=\alpha_{\cdot}^V\star(\Lambda_{\cdot}+ |Y|\star
C_{\cdot}+\frac{1}{2}\varqM_{\cdot})$.
The semimartingale $X^{\Lambda,C}_{\cdot}(Y)=Y_{\cdot} +\Lambda
_{\cdot}+ |Y|\star
C_{\cdot}=Y_{\cdot}+D^{\Lambda,C}_{\cdot}(Y)$ is associated with
the martingale $M_{\cdot}$ and
the finite variation process $-V_{\cdot}^X$ where $V_{\cdot
}^X=V_{\cdot}-D^{\Lambda
,C}_{\cdot}(Y)=(\alpha_{\cdot}^V-1)\star D^{\Lambda,C}_{\cdot
}(Y)+\frac{1}{2}\alpha_{\cdot}^V\star
\varqM_{\cdot}$.
% where $\alpha^V_{\cdot}\in[-1,1]$ is a predictable process such that
%$V=
Since the process $-V_{\cdot}^X+\frac{1}{2}\varqM_{\cdot}=(1-\alpha
_{\cdot}^V)\star
(D^{\Lambda,C}_{\cdot}(Y)+\frac{1}{2}\varqM_{\cdot})$ is an increasing
process, the
semimartingale $X^{\Lambda,C}_{\cdot}(Y)$ is a $\cQ$-submartingale.

\item[(i.b)]Assume now that both processes $e^{{\lo X}_{\cdot}}$ and
$e^{{\lu
X}_{\cdot}}$ are submartingales, where $ {\lo X}_{\cdot} \equiv
X^{\Lambda,C}_{\cdot}(Y)$
and
${\lu X}_{\cdot} \equiv X^{\Lambda,C}_{\cdot}(-Y)$. The processes
${\lo X}_{\cdot}$ and
${\lu
X}_{\cdot}$ satisfy the following relations, where $D^{\Lambda
,C}_{\cdot} \equiv
D^{\Lambda,C}_{\cdot}(Y)$:
\[
\tfrac{1}{2}({\lo X}_{\cdot}-{\lu X}_{\cdot})=Y_{\cdot}
\quad\mbox{and}\quad \tfrac{1}{2}
({\lo X}_{\cdot}+{\lu X}_{\cdot}) = D^{\Lambda,C}_{\cdot}=\Lambda
_{\cdot}+\tfrac{1}{2}|{\lo X}_{\cdot} -
{\lu X}_{\cdot}|\star C_{\cdot}.
\]
Using the same notation and arguments as above, the processes
${\lo X}_{\cdot}$ and ${\lu X}_{\cdot}$, whose exponentials are
submartingales,
can only have positive jumps. This contradicts the fact that their sum
is a continuous increasing process. Hence, both processes are continuous.
For the same reasons, the sum ${\lu M}_{\cdot}+{\lo M}_{\cdot}$ is identically
equal to $0$, and the sum of increasing processes
$\frac{1}{2}({\lu A}_{\cdot}+{\lo A}_{\cdot} )=D^{\Lambda,C}_{\cdot
}+\frac{1}{2}\varq{\lo M}_{\cdot}
\equiv\frac{1}{2}G^{\Lambda,C}_{\cdot}$.

There exists a predictable process $\alpha_{\cdot}$, with $\alpha
_{\cdot}\in
[0,2]$, such that ${\lo A}_{\cdot}=\frac{1}{2}\alpha_{\cdot}\star
G^{\Lambda,C}_{\cdot}$.
Substituting ${\lo A}_{\cdot}$ in the decomposition of $Y_{\cdot
}=\frac{1}{2}({\lo
X}_{\cdot}-{\lu X}_{\cdot})$, we get $dY_{t}=- \frac{1}{2}(1-\alpha
_{t})\,dG^{\Lambda
,C}_t+d{\lo M}_t $.
Therefore, $Y_{\cdot}$ is a $\cQ(\Lambda,C)$-semimartingale.

\item[(ii.a)] Let $Y_{\cdot}$ be a $\cQ(\Lambda, C)$-semimartingale. Since
$X^{\Lambda,C}_{\cdot}(Y)=Y_{\cdot}+D^{\Lambda,C}_{\cdot}$, we have
$e^{Y_{\cdot}}=e^{-D^{\Lambda,C}_{\cdot}}e^{X^{\Lambda,C}_{\cdot
}(Y)}$. From the
classical It\^o's formula,
\[
de^{Y_t}=e^{-D^{\Lambda,C}_t}\,de^{X^{\Lambda
,C}_t(Y)}-e^{Y_t}\,dD^{\Lambda,C}_t\quad
\mbox{and}\quad
dU^{\Lambda,C}_t(e^{Y})=e^{-D^{\Lambda,C}_t}\,de^{X^{\Lambda,C}_t(Y)}.
\]
%}

Then when $Y_{\cdot}$ is a continuous process, $\exp(X^{\Lambda,
C}_{\cdot}(Y))$ is a
submartingale iff $U^{\Lambda,C}_{\cdot}(e^Y)$ is a submartingale.

%C}_{\cdot}(Y))$ and $U^{\Lambda,C}_{\cdot}(e^Y)$ are submartingales.\\
%
\item[(ii.b)] Assume now that both processes $U_{\cdot}(e^Y)$ and $U_{\cdot
}(e^{-Y})$
are l\`adl\`ag submartingales. Let $U_{\cdot}(e^Y)=U_0+\lo{N}_{\cdot
}+\lo{K}_{\cdot}$
and $U_{\cdot}(e^{-Y})=U_0+\lu{N}_{\cdot}+\lu{K}_{\cdot}$ be their
respective additive
decompositions. As before, we can show that the process $Y_{\cdot}$ is
continuous. The previous equivalence yields to the result.\quad\qed
%Since the processes
%$e^{X^{\Lambda, C}_{\cdot}(Y)}=e^{(\Lambda+|Y|\star C)}\star U_{
%$e^{X^{\Lambda, C}_{\cdot}(-Y)}=e^{(\Lambda+|Y|\star C)}\star U_{
%are two
%submartingales, we can apply the previous results to conclude
%that $X^{\Lambda, C}_{\cdot}(\pm Y)$ is a $\cQ(
%and so is $Y_{\cdot}$.
\end{longlist}
\noqed\end{pf*}
%
%%%%%%%%%%%%%%%%%%%%%%%%%%%%%%%%%%%%%%%%%%%%%%%%%%%%%%%%%%
%===========================
%s3 #&#
\section{Exponential uniform integrability and entropic inequalities}\label{subsectioncharactexpoinequality}
%===========================
In the previous section, we have obtained a simple characterization of
${\mathcal
Q}(\Lambda,C)$-semimartingales using an exponential transformation,
leading naturally to positive submartingales defined by their multiplicative
or additive decomposition. Whenever submartingales have good
integrability properties, the existence of an additive decomposition\vadjust{\goodbreak}
is equivalent to the submartingale inequalities. It is the famous
Doob--Meyer decomposition. The main objective of this section is
to precise such integrability properties and the subsequent inequalities.
%%%%%%%%
%s3.1 #&#
\subsection{Uniform integrability, class ($\cD$) and their exponential equivalents}
\label{uniformintegrability} \label{parclassUexp}\label{parclassDexp}
%equivalents}
%%%%%%%%
%pa3.1.0.1 #&#
\subsubsection*{The class ${\mathcal U}_{\exp}$}
In the classical martingale theory, uniformly integrable (u.i.)
martingales (in particular the conditional expectation of some positive
integrable random variable) play a key role as martingale equalities
are then valid between two stopping times.
%Recall that the conditional
%expectation of some positive integrable random variable is still a
%uniformly integrable martingale.
The class of such martingales is denoted by $\cU$.

In the exponential framework, any exponential martingale
$\mathcal{E}(M)_{\cdot}$ of a continuous martingale $M_{\cdot}$ is a positive
local martingale, with expectation $\leq1$, hence a
supermartingale. The process $\mathcal{E}_{\cdot}(M)$ is a u.i. martingale
on $[0,T]$ if and only $\mathcal{E}_t(M) =
\E[\mathcal{E}_T(M)|\mathcal F_t]$ $\mathbb{P} \mbox{-a.s.}$
%Note that, thanks to the supermartingale property, it is enough to
%check the inequality $\mathcal{E}_t(M) \leq\E%[\mathcal{E}_T(M)|
It is therefore natural to introduce the class ${\mathcal U}_{\exp}$
of continuous martingales $M$
such that $\mathcal{E}_{\cdot}(M)$ is a uniformly integrable martingale.
%
%pa3.1.0.2 #&#
\subsubsection*{The classes ${\mathbb L}^1_{\exp}$ and $({\mathcal
D}_{\exp})$}
A $\cF_T$-measurable random variable $X_T$
belongs to ${\mathbb L}^1$ provided that ${\mathbb E}(|X_T|)<\infty$
and by definition belongs to ${\mathbb L}^1_{\exp}$ if $\exp(X_T)\in
{\mathbb L}^1.$

The optional processes $X$ for which the absolute value is
dominated by a uniformly integrable martingale are said to be in the
class\footnote{P. A. Meyer used the term ``class $(\cD)$,'' in the
honor of J. L. Doob.} $(\cD)$. They are also characterized by the fact
that the family of random variables $\{X_\sigma;
\sigma\leq T, \sigma\mbox{ stopping times}\}$ is uniformly
integrable.
When adopting the exponential point of view, we can extend this notion
into:
\[
X_{\cdot}\mbox{ is said to be in the class }(\mathcal
{D}_{\mathrm{exp}})\mbox{ if }e^{X_{\cdot}}\mbox{ belongs to the class }(\mathcal{D}).
\]
Observe that $|X_{\cdot}|$ belongs to the class $(\mathcal{D}_{\mathrm{exp}})$
if and only if $X_{\cdot}$ and $-X_{\cdot}$ belong to the class
$(\mathcal{D}_{\mathrm{exp}})$. The sufficient condition is based on the
intermediate result that $\cosh(X_{\cdot})=\cosh(|X_{\cdot}|)$
is in the class $(\mathcal{D})$.
%
%pa3.1.0.3 #&#
\subsubsection*{$(\mathcal{D})$-submartingales and conditional
inequalities}
%{\cgreen Pauline: c'est plus clair, mais j'ai encore un peu de mal.
%J'ai essaye de modifie un tout petit peu la redaction de la premiere
%phrase pour etre plus claire. A %verifier. En particulier doit-on
%ecrire seulement "sumartingale" ou "$\cQ$-submartingale"?} \\
A submartingale $S_{\cdot}$ (as defined in its general form in
Section \ref
{parrecallsemimartingale}), which is in the
class $(\mathcal{D})$, satisfies the following conditional
``submartingale inequality'' \label{submartingaleinequalities}
\[
\mbox{for any stopping times }\sigma\leq\tau\leq T\qquad
S_\sigma\leq\E[S_\tau|\cF_\sigma],\qquad \mbox{a.s.}
\]

Conversely, it is well known that any c\`{a}dl\`{a}g process in the
class $(\mathcal{D})$ satisfying these inequalities admits a Doob--Meyer
decomposition into a martingale and a predictable c\`{a}dl\`{a}g
increasing process (see Protter~\cite{Protter}, Chapter~3), that is, is
a submartingale in the previous sense. The less standard l\`{a}dl\`{a}g
case, motivated by optimal stopping problems, has been established by
Dellacherie and Meyer~\cite{Dellacherie-Meyer}.\vadjust{\goodbreak}
\subsection{Entropic inequalities and quadratic semimartingales}\label
{subsecentropicinequalities}\label{parentropicq-semimartingale} \label{parentropicq-semimartingale}\label{parmaximalinequalities}\label
{changeprobaentropy}
%%%%%%%%%%%
%{\cgreen Pauline: je trouve cela vraiment tres tres clair maintenant!}
%pa3.2.0.1 #&#
\subsubsection*{Entropic submartingales}
When considering a positive $(\mathcal{D})$-submartin\-gale~$S_{\cdot }$,
the logarithm $X_{\cdot}=\ln{S_{\cdot}} $ is a $\cQ$-submartingale in
the class $(\mathcal{D}_{ \exp})$ and satisfies the so-called
\textit{entropic inequality}:
%e11 #&#
%
\begin{eqnarray}\label{eqentropicinequalities}
\forall\sigma\leq\tau\leq T\qquad X_\sigma\leq\rho_\sigma(X_\tau
)
\nonumber
\\[-8pt]
\\[-8pt]
\eqntext{\mbox{a.s. where }
\rho_{\sigma}(X_\tau)= \ln\mathbb{E}[\exp(X_\tau)|\cF_\sigma].}
\end{eqnarray}
%
%}
The operator $\rho_{\cdot}$ is known as the \textit{entropic
process} and
has been intensively studied in the framework of risk measures (see,
e.g., Barrieu and El Karoui~\cite{Barrieu-ElKaroui5} or
\cite{Barrieu-ElKaroui6}).
Since conversely, any $\cQ$-submartingale in the class $(\mathcal{D}_{
\exp})$ satisfies the entropic inequalities, we refer to it as \textit{entropic submartingale.}

%{\cgreen Pauline: on a parle de ''submartingale in the class $(
%pour aller a la classe $(\mathcal{D}_{
%}\\
An example of entropic submartingale is the simple process $r_{\cdot}(M)$
defined in equation
(\ref{eqrtM}) with $M_{\cdot}\in{\mathcal U}_{\exp}.$
%(Subsection~\ref{classUexp})
In this case,
$\exp r_{\cdot}(M)=\mathcal{E}_{\cdot}(M)$ is a positive u.i.
martingale, equal to
the conditional expectation of its terminal value $\exp(r_T(M))$. Since
$\xi_T \equiv r_T(M) \in{\mathbb L}^1_{\exp}$, we can recover
$r_t(M)$ from its terminal condition from the following
identity\footnote{Note that the identity
$\rho_{t}(\xi_T)=r_t(\rho_{0}(\xi_T),M)$
has suggested the notation $r_t(M)$ for the logarithm of some
exponential martingale.} based on the entropic process $\rho_{\cdot
}(\xi_T)$:
%e12 #&#
%
\begin{equation}\label{eqentropicprocess}
r_t(M)= \ln\mathbb{E}[\exp(\xi_T)|\cF_t]=\rho_t(\xi_T),\qquad
\xi_T \equiv r_T(M).
\end{equation}
The conditional properties of the u.i. martingale $\exp
(r_t(M))=\break\mathbb
{E}[\exp(\xi_T)|\cF_t] =\E[\exp(\xi_T)]\mathcal{E}_t(M)$ are
translated into the time consistency property of the entropic process
over any pair of stopping times $(\sigma,\tau)$ such that
$\sigma\leq\tau$, $\rho_{\sigma}(\xi_T)=\rho_{\sigma}(\rho
_{\tau
}(\xi_T))$.

Finally, let us observe that $\rho_{\cdot}(\xi_T)$ is the smallest
$q$-semimartingale $X_{\cdot}=X_0+r_{\cdot}(N)$ with the terminal
value $ X_T=\xi_T$.
This is a simple consequence of the fact that $\exp(X_{\cdot})$ is a positive
local martingale and hence a supermartingale.
%%%
%pa3.2.0.2 #&#
\subsubsection*{Entropic inequalities and $\cQ$-semimartingales} We are now able to give another
formulation for the characterization of $\cQ$-semimartingales in the
class $(\mathcal{D}_{\exp})$ in terms of inequalities involving the
entropic process.
%, or submartingale inequalities.
This formulation will prove to be better suited than that of Theorem
\ref{thcharacterizationgeneralquasimart} when taking limits as we will
see in a later section. %%
%{\cgreen Je propose de simplifier la suite, en faisant deux theoremes.
%Le premier sur les
%$\cQ$-semimartingales en supprimant les conditions \rmi et \rmii, dont
%on ne se sert pas veritablement.
%Le deuxieme qui deborde sur la suite, et qui donne des conditions
%necessaire et suffisante sur $\bar{X}^{\Lambda,C}_t(|Y|))$ pour etre
%de la classe D}
%th3.1 #&#
%
\begin{theorem}
\label{thcharacterizationqmartDexp} Let $X_{\cdot}$ be
a l\`{a}dl\`{a}g optional process such that $|X_T|\in\mathbb
L^1_{\exp}$. Then $X_{\cdot}$ is a $\cQ$-semimartingale such that
$|X_{\cdot}|\in(\mathcal{D}_{\exp})$ if and only if $X_{\cdot}$
and $-X_{\cdot}$ are
entropic submartingales, or equivalently
if for any pair of stopping
times $0
\leq\sigma\leq\tau\leq T$,
%e13 #&#
%
\begin{equation}\label{eqentropicinequalities}
-\rho_{\sigma}(-X_{\tau}) \equiv{\lu{\rho}}_{\sigma}(X_{\tau
}) \leq X_{\sigma} \leq\rho_{\sigma}(X_{\tau}),\qquad
\mathbb{P}\mathrm{\mbox{-}a.s.}
\end{equation}
\end{theorem}

\begin{pf}
% Thanks to the remarks \rmiv and \rmv of Subsection
Thanks to Section~\ref{uniformintegrability}, when $|X_T|\in
\mathbb
L^1_{\exp}$, the following equivalences hold true [$X_{\cdot}$ is a
$\cQ
$-semimartingale such that
$X_{\cdot}$ and $-X_{\cdot}$ are in the class $(\mathcal{D}_{\exp
})$] is equivalent
to [$e^{X_{\cdot}}$ and $e^{-X_{\cdot}}$ are ($\mathcal
{D}$)-submartingales] that
is equivalent to ($e^{X_{\cdot}}$ and $e^{-X_{\cdot}}$ satisfy the
submartingale
inequalities).
\end{pf}
%
%pa3.2.0.3 #&#
\subsubsection*{Entropic inequalities and $\cQ(\Lambda
,C)$-semimartingales}
The same type of characterization applied to the processes $X^{\Lambda
,C}_T(Y)$ or $U^{\Lambda,C}_T(e^{Y})$ involves inequalities depending
on the process $Y_{\cdot}$ itself and therefore is often difficult to
use. A
possible (but not equivalent) way is to work with the process $\bar
{X}^{\Lambda,C}_{\cdot}(|Y|)$ defined as $\bar{X}^{\Lambda
,C}_t(|Y|) \equiv
e^{C_t}|Y_t|+\int_0^t
e^{C_s}\,d\Lambda_s$ as a generalization of $|Y_{\cdot}|$ by assuming
that the
process $\exp(\bar{X}^{\Lambda,C}_{\cdot}(|Y|))$ is in the class
($\cD$).
%==================================
%pr3.2 #&#
%
\begin{eproposition}
\label{propcharacterizationbarXDexp}
Let $\bar{X}^{\Lambda,C}_t(|Y|) \equiv e^{C_t}|Y_t|+\int_0^t
e^{C_s}\,d\Lambda_s$.
\begin{longlist}[(ii)]
\item[(i)]Let $Y$ be a $\cQ(\Lambda,C)$ semimartingale. Then the process
$\bar{X}^{\Lambda,C}_{\cdot}(|Y|)$ is a $\cQ$-submartingale.

\item[(ii)]Let $Y$ be an optional l\`{a}dl\`{a}g process with $\bar
{X}^{\Lambda,C}_T(|Y|)\in{\mathbb L}^1_{\exp}$.
\end{longlist}

Then the process $\bar{X}^{\Lambda,C}_{\cdot}(|Y|)$
is an entropic submartingale
%$\cQ$-submartingale in the class ($\cD_{\exp}$)
if and only if the following inequalities hold true for any pair of stopping
times $0\leq\sigma\leq\tau\leq T$, where for $t\leq u$, $C_{t,u}
\equiv C_u - C_t$,
%e14 #&#
%
\begin{equation}\label{hypDirectinequality}
|Y_\sigma| \leq\rho_\sigma\biggl(e^{C_{\sigma,\tau}}|Y_\tau|+\int
_\sigma
^\tau e^{C_{\sigma,t}}\,d\Lambda_t\biggr), \qquad\mathbb{P}\mathrm{\mbox{-}a.s.}
\end{equation}
\end{eproposition}

\begin{pf} For the sake of simplicity, we omit $Y$ in $\bar X^{\Lambda
,C}_t(|Y|)$ and $D^{\Lambda,C}_{\cdot}(|Y|)$.

\begin{longlist}[(ii.a)]
\item[(i)]By It\^o--Tanaka formula involving the sign function [$\operatorname{sign}(x)=x/|x|$], with $\operatorname{sign}(0)=0$, and the local time $ L_{\cdot}(Y)$
of $Y.$ at $0$, $|Y_{\cdot}|=|Y_0|+\operatorname{sign}(Y)\star
Y_{\cdot}+L_{\cdot}(Y)=|Y_0|+M^s_{\cdot}-V^s_{\cdot}+L_{\cdot
}(Y)$, where $dM^s_t=\operatorname{sign}(Y)_t
\,dM_t$ and $dV^s_t=\operatorname{sign}(Y)_t\, dV_t$.
This decomposition leads to the following representation of the
differential of $\bar X^{\Lambda,C}_{\cdot}= e^{C_{\cdot}} |Y_{\cdot
}|+e^{C_{\cdot}}\star
\Lambda_{\cdot}$:
\begin{eqnarray*}
d\bar X^{\Lambda,C}_t&=&e^{C_t}[|Y_t|\,dC_t+d\Lambda_t+dM^s_t-
dV^s_t+dL_t(Y)]\\
&=&e^{C_t}\bigl[dD^{\Lambda,C}_t+\tfrac{1}{2}\,d\varqM_t-dV^s_t+dL_t(Y)
\bigr]+e^{C_t}\bigl(dM^s_t-\tfrac{1}{2}\,d\varqM_t\bigr).
\end{eqnarray*}

Observe that $\bar A^s_{\cdot}=D^{\Lambda,C}_{\cdot}+\frac
{1}{2}\varqM
_{\cdot}-V^s_{\cdot}+L_{\cdot}(Y)$
is an increasing process. The martingale part of $\bar X^{\Lambda
,C}_{\cdot}(|Y|)$ is $\bar M^{C}_{\cdot}=e^{C_{\cdot}}\star
M^s_{\cdot}$ with quadratic
variation $\varq{\bar M^{C}}_{\cdot}=e^{2C_{\cdot}}\star\varqM
_{\cdot}$. So, the
following decomposition shows that $\bar X^{\Lambda,C}_{\cdot}$ is a
$\cQ
$-submartingale since $e^{C_{\cdot}}-1\geq0 $,
\[
d\bar X^{\Lambda,C}_{\cdot}=e^{C_{\cdot}}\bigl[d\bar
A^s_{\cdot}+\tfrac{1}{2}
(e^{C_{\cdot}}-1)\,d\varqM_{\cdot}\bigr] + d r_{\cdot}(e^{C_{\cdot
}}\star
M^s_{\cdot}).
\]

\item[(ii.a)] The assumption that $\exp(\bar{X}^{\Lambda,C}_{\cdot})$ is
a ($\cD
$)-submartingale
implies in particular that
$\bar{X}^{\Lambda,C}_T\in{\mathbb L}^1_{\exp}$, and that
$\bar{X}^{\Lambda,C}_0=|Y_0|\leq\rho_0(\bar{X}^{\Lambda,C}_T)$. The
same inequality
holds true if\vadjust{\goodbreak} we start at time $\sigma$ with horizon $\tau$ by
considering the $\sigma$-conditional
expectation of $\bar{X}^{\Lambda,C}_{\sigma,\tau}=e^{C_{\sigma
,\tau
}}|Y_\tau|+\int_\sigma^\tau e^{C_{\sigma,t}}\,d\Lambda_t$, so that
$|Y_\sigma| \leq\rho_\sigma(\bar{X}^{\Lambda,C}_{\sigma,\tau}).$

\item[(ii.b)] Conversely, assume inequality (\ref{hypDirectinequality}),
$|Y_\sigma
| \leq\rho_\sigma(\bar{X}^{\Lambda,C}_{\sigma,\tau}).$
Observe that the entropic process
$\rho_{\delta,t}(\xi_T)=\frac{1}{\delta}\rho_{t}(\delta\xi_T)$ is
increasing with respect to the parameter $\delta$ (from the H\"{o}lder
inequality for the exponential). Then, since $e^{C_\sigma}\geq1$, we have:
$\rho_\sigma(e^{C_{\sigma,\tau}}|Y_\tau|+\int_\sigma^\tau
e^{C_{\sigma
,t}}\,d\Lambda_t)\leq e^{-C_\sigma}
\rho_\sigma(e^{C_{0,\tau}}|Y_\tau|+\int_\sigma^\tau
e^{C_{0,t}}\,d\Lambda_t)
$. So $\bar{X}^{\Lambda,C}_{\cdot}=e^{C_{\cdot}} |Y_{\cdot}|+
e^{C}\star\Lambda_{\cdot}$
satisfies the entropic inequalities $\bar{X}^{\Lambda,C}_\sigma
\leq\rho_\sigma(e^{C_{0,\tau}}|Y_\tau|+\int_0^\tau
e^{C_{0,t}}\,d\Lambda_t+\int_0^\tau e^{C_{0,t}}\,d\Lambda_t)= \rho
_\sigma(\bar{X}^{\Lambda,C}_\tau)$. Taking $\tau=T$, it follows that
$\bar{X}^{\Lambda,C}_{\cdot}$ is dominated by the $(\cD_{\exp
})$-process $\rho
_{\cdot}(\bar{X}^{\Lambda,C}_T)$ and so is an entropic submartingale.
%dominated on $[0,\tau]$ by the entropic process
%$\rho_{\cdot}(\bar{X}^{\Lambda,C}_\tau(|Y|))$ in $(\cD_{\exp})$ since $
Hence, the result.\quad\qed
\end{longlist}
\noqed\end{pf}

The properties of the dominating process $\rho_{\cdot
}(e^{C_{\cdot,T}}|Y_T|+\int
_{\cdot}^Te^{C_{\cdot,s}}\,d\Lambda_s)$ are therefore essential to
obtain results for the process $Y_{\cdot}$. The nonadapted process
$\phi_{\cdot,T}(|Y_T|)=e^{C_{\cdot,T}}|Y_T|+\int_{\cdot
}^Te^{C_{\cdot,s}}\,d\Lambda_s$ with
initial condition $\phi_{0,T}(|Y_T|)= \bar X^{\Lambda,C}_{T}(|Y|)$,
first introduced in Briand and Hu~\cite{Briand-Hu}, Lemma 1, is the
positive decreasing solution
of the ordinary differential equation, with terminal condition
$|Y_T|$,
%e15 #&#
%
\begin{equation}\label{eqphi}
d\phi_t=
-(d\Lambda_t+|\phi_t|\,dC_t),\qquad \phi_T=|Y_T|.
\end{equation}
In order words, the nonadapted process $U^{\Lambda,C}_{\cdot
}(e^{\phi_{\cdot,T}
})=e^{\phi_{\cdot,T} }+\int_0^{\cdot}e^{\phi_{s,T}}\,d\Lambda_s+\int_0^{\cdot}
e^{\phi
_{s,T}}|\phi_{s,T}|\,dC_s$
is constant and equal to $e^{\phi_{0,T} }$. This property is the main
motivation for introducing the $U^{\Lambda,C}$ transformation.

The decreasing property of $\exp(\phi_{\cdot,T})$ explains the
supermartingale property of the process $\Phi_{\cdot}(|Y_T|)$ defined
as the
optional projection of $\exp(\phi_{\cdot,T})$:
%e16 #&#
%
\begin{eqnarray}\label{eqPhi}
\Phi_\sigma(|Y_T|)&\equiv&\E[\exp(\phi_{\sigma,T} (|Y_T|))|\cF
_\sigma]
\nonumber
\\[-8pt]
\\[-8pt]
\nonumber
&=&
\exp\biggl(\rho_\sigma\biggl(e^{C_{\sigma,T}}|Y_T|+\int_\sigma^Te^{C_{\sigma
,t}}\,d\Lambda_t\biggr)\biggr).
\end{eqnarray}
Note that, for the sake of clarity, we often omit the reference to
$Y_T$ in
$\phi_{\cdot,T} (|Y_T|)$, $\Phi_{\cdot}(|Y_T|)$ or $\bar X^{\Lambda
,C}_{T}(|Y_T|)$.
%
%th3.3 #&#
%
\begin{theorem}\label{propUtransform}
Assume $\E[\exp(\bar X^{\Lambda,C}_{T}(|Y_T|))]=\E[\exp(\phi
_{0,T})]<\infty$.

\begin{longlist}[(iii)]
\item[(i)]The process $\Phi_{\cdot}$ is a $(\cD)$-supermartingale
dominated by the
martingale $ \E[e^{\phi_{0,T}}|\cF_t]=N^0_t$, with the additive
decomposition $\Phi_{\cdot}=\Phi_0+N^\Phi_{\cdot} -A^\Phi_{\cdot
}$. The predictable
increasing process is
$A^\Phi_{\cdot} =\int_0^{\cdot} \Phi_s \,d\Lambda_s+\int_0^{\cdot} \E[e^{\phi
_{s,T}}|\phi
_{s,T}||\cF_s]\,dC_s$, when the process $N^\Phi_{\cdot}$
is a uniformly integrable martingale.

\item[(ii)]The process $U^{\Lambda,C}_{\cdot}(\Phi)=\Phi_{\cdot}+\int
_0^{\cdot}\Phi_s
\,d\Lambda_s+ \int_0^{\cdot}\Phi_s\ln(\Phi_s)\,dC_s$ is a positive
$(\cD)$-supermartingale, associated with the same u.i. martingale
$N^\Phi_{\cdot}$, and the increasing process $A^U_{\cdot}=\int_0^{\cdot}
(\E[e^{\phi_{s,T}}|\phi_{s,T}||\cF_s]-\Phi_s\ln(\Phi_s))\,dC_s$.

\item[(iii)]Assume inequality (\ref{hypDirectinequality}) for the process
$|Y_{\cdot}|$. The
processes $U^{\Lambda,C}_{\cdot}(e^Y)$ and $U^{\Lambda,C}_{\cdot
}(e^{-Y})$ are two
$(\mathcal{D})$-submar\-tingales dominated by the $(\mathcal
{D})$-supermar\-tingale
$U^{\Lambda,C}_{\cdot}(\Phi)$.\vspace*{-1pt}
\end{longlist}
\end{theorem}

%re1 #&#
%
\begin{remark}\label{remUtransform}
The positive quantity $H^{\mathrm{ent}}_s(e^{\phi_{s,T}}) \equiv
\E[e^{\phi_{s,T}}\phi_{s,T}|\cF_s]-\break\Phi_s\ln(\Phi_s)$ appearing in
$A^U_{\cdot}$ is well known in statistics as the conditional Shannon entropy
of the random variable $e^{\phi_{s,T}}$. Its properties will be studied
in the next subsection when considering integrability properties of the
supremum.\vspace*{-1pt}
\end{remark}

\begin{pf*}{Proof of Theorem~\ref{propUtransform}}
As observed by Briand and Hu~\cite{Briand-Hu}, Lemma~1, since $\phi_{t,T}$
is a positive solution of the differential equation $d\phi_t=
-(d\Lambda_t+\phi_t \,dC_t)$, the nonadapted process $U^{\Lambda
,C}_t(e^{\phi_{\cdot,T} })$ is constant,
$U_t^{\Lambda, C}(e^{\phi_{\cdot,T}})=e^{\phi_{t,T} }+\int_0^te^{\phi
_{s,T}}\,d\Lambda_s+\int_0^t e^{\phi_{s,T}}|\phi_{s,T}|\,
dC_s=e^{\phi_{0,T}},$
with $\phi_{0,T} =\bar{X}^{\Lambda,C}_T\in\mathbb L^1_{\exp}$.
The dynamics of the supermartingale
$\Phi_t=\E[e^{\phi_{t,T} }|\cF_t]$ is obtained by taking conditional
expectation in this relation.\vspace*{-1pt}

\begin{longlist}[(iii)]
\item[(i)]First, observe that the assumption $e^{\phi_{0,T}}\in{\mathbb
L}^1$ implies that $e^{\phi_{T,T}}\in{\mathbb L}^1$ and
that the nonadapted increasing process $B^\phi_t=\int_0^te^{\phi
_{s,T}}\,d\Lambda_s+\break\int_0^t e^{\phi_{s,T}}|\phi_{s,T}|\,dC_s$ is
integrable.
Since $\Phi_{\cdot}$ is the optional projection of $e^{\phi_{\cdot,T}}$, and
since both increasing processes $\Lambda_{\cdot}$ and $C_{\cdot}$
are adapted, the
dual predictable projection of $B^\phi_t$ is the continuous process
$A^\Phi_t=\int_0^t \Phi_s \,d\Lambda_s+\break\int_0^t\E[e^{\phi
_{s,T}}|\phi
_{s,T}||\cF_s]\,dC_s$, generating the same conditional variation, $\E
[B^\phi_{t,T}-A^\Phi_{t,T}|\cF_t]=0 $.
%
%{\cgblue that is?} by definition has the same
%expectation than $B^\phi_T$, and also belongs to $ {\mathbb L}^1$.
So the
process $N^1_t=\E[B^\phi_T-A^\Phi_T|\cF_t]=\E[B^\phi_t-A^\Phi
_t|\cF_t]$ is
a uniformly integrable martingale. Then, taking the conditional
expectation of the constant process $U_{\cdot}^{C,\Lambda}(e^{\phi_{\cdot,T}})$
implies that
$\Phi_t+A^\Phi_t+N^1_t=N^0_t$, and
$N^\Phi_t=N^0_t-N^1_t$.

\item[(ii)]To show that $U^{\Lambda,C}_{\cdot}(\Phi)$ is also a
supermartingale, we use that the Shannon entropy (see Remark \ref
{propUtransform}) $H^{\mathrm{ent}}_s(e^{\phi_{s,T}})=\E[e^{\phi
_{s,T}}|\phi
_{s,T}||\cF_s]-\Phi_s\ln(\Phi_s)$ is positive, and the process
$A^U_{\cdot}=\int_0^{\cdot}H^{\mathrm{ent}}_s(e^{\phi_{s,T}})\,dC_s$ is
increasing. Then,
some simple
calculation shows that $U^{\Lambda,C}_{\cdot}(\Phi)+A^U_{\cdot
}=\Phi_{\cdot}+A^\Phi
_{\cdot}=N^\Phi_{\cdot}
$ is a positive u.i. martingale, that provides the Doob--Meyer
decomposition of the
supermartingale $U^{\Lambda,C}_{\cdot}(\Phi)$.

\item[(iii)]This last statement is a straightforward consequence
of the inequality $e^{|Y_{\cdot}|}\leq\Phi_{\cdot}$.\quad\qed\vspace*{-1pt}
\end{longlist}
\noqed\end{pf*}

%%%
%re2 #&#
%
\begin{remark}\label{remetaT}
The key condition to obtain these properties is that the process
$U^{\Lambda,C}_{\cdot}(\Phi(|Y_T|))$ is a ($\cD$)-supermartingale.
Note that this is also true if we replace $|Y_T|$ by
any $\cF_T$-random variable $|\eta_T|\geq|Y_T|$, such that
$e^{C_{T}}|\eta_T|+\int_0^Te^{C_{s}}\,d\Lambda_s \in{\mathbb
L}^1_{\exp}$.\vspace*{-1pt}
%Therefore, in the next
%section, we will consider a slightly modified version of Hypothesis
%random variable $\eta_T$.
\end{remark}

%re3 #&#
%
\begin{remark} As observed by Briand and Hu~\cite{Briand-Hu}, extending
the results of Lepeltier and San Martin in~\cite{LepeltierSanMartin98},
the linear growth condition in\vadjust{\goodbreak} $Y_{\cdot}$, $|Y_{\cdot}|\star
C_{\cdot}$,
may replaced by a superlinear growth $h(|Y_{\cdot}|)\star C_{\cdot}$,
where $h$ is
an increasing convex $C^1$ function, with $h(0)>0$,
satisfying the integrability condition $\int_0^T du \frac
{|u|}{h(u)}=+\infty$.
The function $\phi(t)$ is then replaced by the solution of the ODE
$\phi'(t)=- h(\phi_t)$ with a terminal condition $\phi(T)=z \geq0$.
\end{remark}

%======================
%pa3.2.0.4 #&#
\subsubsection*{Maximal exponential integrability and $L\log
L$-condition} When\footnote{This
paragraph can be omitted for a first reading.} looking for entropic
inequalities, assuming that the exponential of the processes is in the
class~($\mathcal{D}$) is a minimal assumption. However, it is sometimes
interesting to obtain estimates on the exponential of the maximum of
these processes. Entropic inequalities reduce the problem to the
estimation of the running supremum of some entropic processes, or
equivalently to the running supremum of some positive martingales, for
which we can apply standard Burkholder--Davis--Gundy (BDG) martingale
inequalities. An excellent presentation of the different martingale
inequalities may be found in Lenglart, Lepingle and Pratelli
\cite{Lenglart-Lepingle-Pratelli}.

From now on, we adopt the following nonstandard notation for the
running supremum of some measurable process $X$: $\max|X_{t}|=\max
_{0\leq u\leq t}|X_u|$ and $\max|X_{s,t}|=\max_{s\leq u\leq
t}|X_u-X_s|$. The space of semimartingales $X_{\cdot}$ such that $\max
|X_{T}|\in\mathbb L^p$ $(p\geq1)$ is denoted by $\mathcal S^p$. For
continuous local martingales, the relevant quantity is the quadratic
variation and we denote by $\mathbb{H}^p$ the space of martingales with
a quadratic variation in $\mathbb{L}^p$. Moreover, for any continuous
local martingale $M_{\cdot}$, such that $M_0=0$, the BDG inequality gives
some estimates of its maximum in terms of with its quadratic variation
as, for any $0<p<\infty$, there exist two positive constants $c_p$ and
$C_p$ such that:
\[
\mbox{for any }0<p<\infty\qquad  c_p \mathbb E[\varqM_T^{p/2}]\leq
\mathbb E[\max|M|_{T}^{p}]\leq C_p \mathbb E[\varqM_T^{p/2}].
\]

The following Doob inequalities, based on the terminal condition and
only true for $p>1$, are more classical:
\[
\mbox{for any }p>1\qquad k_p \mathbb E[\varqM_T^{p/2}]\leq
\mathbb
E[ |M|_{T}^{p}]\leq K_p \mathbb E[\varqM_T^{p/2}].
\]

So, for $p>1$, $\mathbb E[\max|M|_{T}^{p}]<\infty$ if and only
$\mathbb E[\varqM_T^{p/2}]<\infty$. In other words, the spaces
$\mathcal S^p$ and
$\mathbb H^p$ coincide.

In terms of exponential martingale $L_{\cdot}=\mathcal E_{\cdot}(M)$, these
results become
\begin{eqnarray*}
%when $p>1$:
\forall p>1 \qquad L_{\cdot} \in\mathcal S^p
\quad&\Longleftrightarrow&\quad L_T=\exp(r_T(M))\in\mathbb L^p\\
\quad&\Longleftrightarrow&\quad
\biggl(\int_0^TL^2_s\,d\varqM_s\biggr)^{1/2}\in\mathbb L^p.
\end{eqnarray*}

When $p<1$, a similar maximal inequality holds true for exponential
martingales or more generally for positive supermartingales \cite
{Lenglart-Lepingle-Pratelli},
%such that $L_0=1$ as
\[
\forall p<1\qquad \mathbb E[\max L_{T}^{p}]\leq\frac
{\mathbb E((L_0)^p)}{1-p}.\vadjust{\goodbreak}
\]

When $p=1$ and the local martingale is positive, we have to use the
following $L\log L$-condition.
%that we recall in following form due to Harremo\"es~\cite{Harremoes}
%in the discrete time.
%The so-called $L\log L$-Doob inequality (see e.g. Protter
%gives
%a necessary and sufficient condition on the terminal value of a u.i.
%exponential martingale $L_{\cdot}={\mathcal E}(M)_{\cdot}$
%for it to be in ${\mathbb H}^1$, or equivalently to have an integrable
%running supremum
%$\max L_t\equiv\sup_{0\leq s\leq t}L_s$.
% More precise sharp estimates have been recently
%proposed by Harremo\"es~\cite{Harremoes} in the discrete time
%context.We give a similar result in the continuous case.
%
%pr3.4 #&#
%
\begin{eproposition} \label{propLlogLinequality} Let $L_{\cdot}=\exp
(M_{\cdot}-\frac{1}{2}\varqM_{\cdot})$ be a positive continuous locale
martingale and
$\max L_t$ its running supremum.

%For any $m>0$, let $u_m(x)$ be the convex function defined on $
%%m$ is attained at $x=m$. The classical case corresponds to $m=1$, and
%then the function $u_1$ is just denoted by $u$.\\
\begin{longlist}[(iii)]
\item[(i)](Doob) Assume that $L_{\cdot}$ is a u.i. martingale.
Then
\[
\mathbb E(\max L_T)-1=\mathbb E(L_T\ln(\max
L_T))\geq\mathbb E(L_T \ln(L_T))
\]
and
\[\mathbb E(L_T \ln(L_T))=\mathbb E\bigl(L_T
\tfrac{1}{2}\varqM_T\bigr).
\]

%
%with $L_0=1$.
% Then $\max L_T$ is an integrable variable if and only if $H^{\rm
%ent}(L_T)=\mathbb E(L_T \ln(L_T))<\infty$. Moreover $H^{\rm
%ent}(L_T)=\mathbb E(L_T \ln(L_T))=\mathbb E(L_T \demi
% \rmb Write $L_t=\mathbb E(L_T|\mathcal F_t)$ in the multificative
%form $L_t=\exp(M_t-\demi\varqM_t)$
% Using the {\cgblue multiplicative decomposition} of the martingale
%$L_t=\mathbb E(L_T|\mathcal F_t)$ as $L_t=\exp(M_t-\demi\varqM_t)$,
% \begin{equation}\label{eqentropyvariance}
% H^{\mathrm{ent}}(L_T)=\mathbb E(L_T \ln(L_T))=\mathbb E(L_T
%
\item[(ii)](Harremo\"es) The following inequality is sharp:
%e17 #&#
%
\begin{equation}\label{eqLlogLinequality}
\mathbb E(\max L_T)-1-\ln(\mathbb E(\max L_T))\leq\mathbb
E(L_T \ln(L_T))=H^{\mathrm{ent}}(L_T).
\end{equation}
The martingale $L_{\cdot}$ belongs to $\mathcal S^1$ if and only if
$\mathbb
E(L_T \ln(L_T))< \infty$.

\item[(iii)]Let $U_{\cdot}$ be a positive $(\cD)$-submartingale with deterministic
initial condition $U_0$ and $m= \mathbb E(U_T)\geq U_0$. The previous
Harremo\"es inequality becomes, when $u_m(x)=x-m-m \ln(x),$
\[
u_m({\mathbb E}(\max U_T))-u_m(U_0)\leq{\mathbb
E}(U_T \ln(U_T))- {\mathbb E}(U_T)\ln(\mathbb E(U_T)
)=H^{\mathrm{ent}}(U_T).
\]

In particular, ${\mathbb E}(\max U_T)$ is dominated by an increasing
function of\break $H^{\mathrm{ent}}(U_T)+u_m(U_0)$.
\end{longlist}
\end{eproposition}

\begin{pf}
The proof is based on Dellacherie~\cite{Dellacherie} and Harremo\"es
\cite{Harremoes}.
\begin{longlist}[(ii.b)]
\item[(i)]Since $L_{\cdot}$ is a continuous process, $\max L_{\cdot}$ only
increases
on the set $\{ L_{\cdot}=\max L_{\cdot}\}$ and $\max L_t=1+\int_0^t
d\max
L_s=1+\int_0^t \frac{L_s}{\max L_s}\,d\max L_s$. Taking the
expectation (after stopping at some stopping time bounding $\max
L_{\cdot}$
on $[0,T]$ if necessary) and using the fact that $L_{\cdot}$ is the
conditional expectation of its terminal value leads to $ \mathbb E(
\max L_T)-1=\mathbb E(L_T\ln(\max L_T))$.

\item[(i.a)] Since $\ln(\max L_T)\!\geq\!\ln( L_T)^+$, and $L_t \ln
(L_t)^-\!\leq\!
1/e$, then $|L_T\ln(L_T)|\!\in\break{\mathbb L}^1$
when $\max L_T\in{\mathbb L}^1$. This establishes the necessary
condition.

\item[(ii)]To prove that finite entropy implies integrability of the $\max$,
we show inequality \eqref{eqLlogLinequality}. We start by studying
$\mathbb E(L_T\ln(\max L_T))-\mathbb
E(L_T\ln( L_T))$ from the concavity of the function $\ln$.
Given that $x^* =\mathbb E(\max L_T)=\break\mathbb
E_{\mathbb Q}(\max L_T/L_T)$ if $\mathbb
Q=L_T.\mathbb P$, $\mathbb E_{\mathbb Q}( \ln(\max L_T/L_T))\leq
\ln(\mathbb E_{\mathbb Q}(\max L_T/L_T))=\ln x^* $.
Inequality \eqref{eqLlogLinequality} is then easily obtained.
An example of c\`{a}dl\`{a}g martingale satisfying the equality may be
found in Harremo\"es~\cite{Harremoes}.

\item[(iii)]The extension to $U_{\cdot}$ being a positive submartingale does
not present any specific difficulties other than purely
computational, since $ \mathbb E( \max U_T)-U_0\leq\mathbb
E(U_T\ln(\max U_T/U_0))$. Taking now $\mathbb
Q=(U_T/m).\mathbb P$, $x^* /m=\break\mathbb E(\max U_T)/m=\mathbb
E_{\mathbb Q}(\max L_T/L_T)$, the convexity inequality becomes:\break
$\mathbb E_{\mathbb Q}(\ln(\max U_T/U_T))\leq\ln(\mathbb E_{\mathbb
Q}(\max U_T/U_T))=\ln(x^* /m)$.\vadjust{\goodbreak} Some elementary algebra gives
the final result.
Observe that $u_m$ is convex and minimal at $z=m$. Since $m_0\leq m$,
$u_m(U_0)\geq u_m(m)$. Then since the entropy is positive, $u_m(U_0)+
H^{\mathrm{ent}}(U_T)$ belongs to the range of $\{u_m(z); z\geq m\}$ and $
\mathbb E( \max U_T)\leq u_m^{-1}(u_m(U_0)+ H^{\mathrm{ent}}(U_T))$.

\item[(i.b)] We now show the link between entropy and quadratic variation.
Assume that $L_T\ln(L_T)\in{\mathbb L}^1$. Let $T_K$ be an increasing
sequence of stop\-ping~times, such that $\ln(L_t)=M_t-\frac{1}{2}\varqM
_t$ is
bounded by $K$. The sequence $T_K$ is increasing and goes to infinity
with $K$. Thanks to the Girsanov theorem, $N^{\mathbb Q}_{\cdot
}=M_{\cdot}-\varqM
_{\cdot}$ is a local martingale with respect to the probability measure
$\mathbb Q=L_T.\mathbb P$, and
$\mathbb E(L_T\frac{1}{2}\varqM_T)=\lim_K \mathbb E(L_{T}\frac{1}{2}
\varqM_{T\wedge T_K})=\break\lim_K\mathbb E(L_{T\wedge T_K}\frac{1}{2}
\varqM_{T\wedge T_K}).$ Using \mbox{$\mathbb E(L_{T\wedge
T_K}N^{\mathbb Q}_{T\wedge T_K})=0$,}
%=\mathbb E(L_{T\wedge T_K}(M_{T\wedge T_K}-\varqM_{T\wedge T_K})
%)=0$,}
%
\begin{eqnarray*}
\mathbb E\bigl(L_{T\wedge T_K}\tfrac{1}{2}\varqM_{T\wedge T_K}\bigr)&=&\mathbb
E\bigl(L_{T\wedge T_K}\bigl(M_{T\wedge T_K}-\varqM_{T\wedge T_K}+\tfrac
{1}{2}\varqM
_{T\wedge T_K}\bigr)\bigr)\\
&=&\mathbb E(L_{T\wedge T_K}\ln(L_{T\wedge T_K}))\leq\mathbb
E( \max L_{T\wedge T_K})-1\\
&\leq&
\mathbb E( \max L_T)-1.
\end{eqnarray*}
Then $N^{\mathbb Q}_{\cdot}$ is a square integrable $\mathbb Q$-martingale
and ${\mathbb E}_{\mathbb Q}(\ln(L_{T}))={\mathbb E}_{\mathbb Q}(\frac{1}{2}
\varqM_T)$, which is is the desired equality.\quad\qed
\end{longlist}
\noqed\end{pf}

Let us now come back to the question of maximal inequalities for
$\mathcal{Q}(\Lambda,C)$-semimartingales.
%such that $\bar{X}^{\Lambda,C}_T(|Y_T|) \in{\mathbb L}^1_{\exp}$.
The various results are based on the behaviour of the entropic process
$\rho_{\cdot}(\bar{X}^{\Lambda,C}_T(|Y_T|))$ also denoted $\rho
_{\cdot}(\bar
{X}^{\Lambda,C}_T).$
% for the sake of notational simplicity.
To give a concise form to the various but similar estimates, we
introduce the following family of positive increasing functions $\psi
_p$ defined on $\mathbb R^+$ by $\psi_p(z)=z^p$ if $p\neq1$ and $\psi
_1(z)=z \ln z-z+1$. Note that, as in the previous subsections, we
consider separately the case of entropic submartingales.
%
%pr3.5 #&#
%
\begin{eproposition} \textup{(i)} Assume $X_{\cdot}=X_0-V_{\cdot}+M_{\cdot}$
to be an entropic
submartingale ($|X_T|\in\mathbb L^1_{\exp}$), such that $\psi_p(
\exp
X_T)\in\mathbb L^1$ provided that $p\geq1$.
%where $\phi_p$ are convex functions defined on $\mathbb R^+$ by$
Then, both processes $\exp(X_{\cdot})$ and ${\mathcal E}(M)_{\cdot}$
belong to
$\mathcal S^p$, and their $\mathcal S^p$ norm are dominated by some
increasing function of $\mathbb E(\psi_p( \exp X_T))$ for $p\geq1$,
and of $\psi_p( \mathbb E(\exp X_T))$ for $p<1$.\vspace*{-6pt}

\begin{longlist}
\item[(ii)]Let $Y_{\cdot}$ be a $\mathcal{Q}(\Lambda,C)$-semimartingale
such that
$\psi_p( \exp\bar{X}^{\Lambda,C}_T)\in\mathbb L^1$ when $p\geq1$.
The processes $\exp(\rho_{\cdot}(\bar{X}^{\Lambda,C}_T))$, $\Phi
_{\cdot}(|Y_T|)$,
$\exp(e^{C_{\cdot}}|Y_{\cdot}|+\int_0^{\cdot} e^{C_s}\,d\Lambda_s)$ and
${\mathcal
E}(e^C\ast M)_{\cdot}$ belong to $\mathcal S^p$ and their $\mathcal S^p$
norms are dominated by some increasing function of $\mathbb E(\psi_p(
\exp\bar{X}^{\Lambda,C}_T))$ for $p \geq1$ or $\psi_p( \mathbb
E(\exp
\bar{X}^{\Lambda,C}_T))$ for $p<1$.
\end{longlist}
\end{eproposition}

\begin{pf} (i) The proof relies on the multiplicative decomposition of
the submartingale $\exp(X_{\cdot})=\exp(X_0+A_{\cdot}){\mathcal
E}(M)_{\cdot}$. Then $\exp
(X_{\cdot})$ and ${\mathcal E}(M)_{\cdot}$ have the same maximal properties.
The proof is a simple consequence of the entropic inequalities (\ref
{propcharacterizationbarXDexp}), BDG inequalities and the maximal
estimates given in Proposition~\ref{propLlogLinequality};\vadjust{\goodbreak}

\begin{longlist}[(ii)]
\item[(ii)]The maximal estimates of $\exp(\rho_{\cdot}(\bar{X}^{\Lambda,C}_T))$
are a simple consequence of (i), and yield to the other estimates
since the different process are dominated by $\exp(\rho_{\cdot}(\bar
{X}^{\Lambda,C}_T))$. For the process ${\mathcal E}(e^C\star M)_{\cdot
}$, we
have to use the decomposition of the entropic submartingale
$e^{C_{\cdot}}Y_{\cdot}+\int_0^{\cdot}e^{C_s}\,d\Lambda_s$.\quad\qed
\end{longlist}
\noqed\end{pf}
%

%pa3.2.0.5 #&#
\subsubsection*{Change of probability measures and entropy} Let $L$ be a positive local martingale with
$L_0=1$. The condition $\mathbb E(L_T \ln(L_T))=H^{\mathrm{ent}}(L)<
\infty$
% is related to the Shannon entropy introduced in Remark
naturally appears when considering the martingale $L$ as
the likelihood of a probability measure $\mathbb{Q}$ equivalent to
$\mathbb{P}$,
%continuous with respect to $\mathbb{P}$,
as it measures the positive
Shannon entropy $H^{\mathrm{ent}}(d\mathbb{Q}/d\mathbb{P})=\mathbb E
(d\mathbb{Q}/d\mathbb{P} \ln(d\mathbb{Q}/d\mathbb{P}))$ of
$\mathbb
{Q}$ with respect to $\mathbb{P}$. The previous result states that
$H^{\mathrm{ent}}(d\mathbb{Q}/d\mathbb{P})$ is finite if and only if the
martingale density $L_{\cdot}$ is in $\mathcal S^1$.

This interpretation is particularly interesting when using the
variational formulation of the the entropic risk measure $\rho_0(\xi
_T)$ (see, e.g., Frittelli~\cite{Frittelli}, F\"ollmer and
Schied~\cite{Foellmer-Schied})
%e18 #&#
%
\begin{equation}\label{eqdualentropy}
\rho_0(\xi_T)=\sup_{\mathbb{Q}}\{\mathbb E_{\mathbb{Q}}(\xi
_T)-H^{\mathrm{ent}}({\mathbb{Q}}/{\mathbb{P}}) | H({\mathbb{Q}}/{\mathbb
{P}})<+\infty
\}.
\end{equation}
In other words, when $\xi_T \in\mathbb L^1_{\exp}$, for any martingale
density $L^Q \in\mathcal S^1$ whose the $\mathcal S^1$-norm is bounded
by $K$, we have an uniform estimate of $\mathbb E_Q(\xi_T)$ given by
$\mathbb E_Q(\xi_T)\leq\rho_0(\xi_T)+\mathbb E(L_T\ln(L_T))\leq
\rho
_0(\xi_T)+K$.

Moreover, when the random variable $\xi_T$ itself is associated with
a finite relative entropy probability measure $\mathbb{Q}^{\xi_T}$
defined by its density $L_T^{\xi_T}=e^{(\xi_T-\rho_0(\xi_T))}$, we
can prove by a simple verification that the supremum is attained for
$\mathbb{Q}^{\xi_T}$.
%In the case when $\mathbb{Q}^{\xi_T}$ does not
%have a finite relative entropy, but $\xi_T$ is \textit{bounded from
%below}, we can approximate $\rho_0(\xi_T)$ by the increasing
%sequence $\rho_0(\xi_T\wedge n)$, and prove that the supremum in the
%optimization problem (\ref{eqdualentropy}) may be restricted to the
%family of \textit{equivalent} probability measures. Such an assumption
%appears in different papers
%dealing with the question of pricing in incomplete markets and entropy.
%(see e.g. Delbaen, Grandits, Rheinlander, Samperi, Schweizer
%and Stricker~\cite{DGRSSS} or Frittelli~\cite{Frittelli}).
Very recently, Choulli and Schweizer~\cite{ChoulliSchweizer} have
developed applications to mathematical finance of the $L\log L$ condition.
%
%NEK: Le commentaire et la proposition qui suivent sont a revoir compte
%tenu de la suite (toujours a faire).\\ Pauline: bien faire le lien
%avec les $\mathcal{Q}(\Lambda,C)$-semimartingales pour la transition
%entre les sections 3 et 4.}
%The following proposition, coming as a straightforward
%consequence of Proposition~\ref{propLlogLinequality}, gives a
%potential application of the $L \log L$ inequality to our framework
%and enables us to have some conditions for the integrability of the
%running maximum process of different processes of interest.
%Then, all the quantity of interest are bounded in $\mathbb L^1$. $
%%%
%The processes $\exp(\rho_{\cdot}(\bar{X}^{\Lambda,C}_T))$, $\Phi_{
%S^p.$\\[0mm]
%e^{C_s}\,d
%are dominated by the Shannon entropy of $\exp(\bar{X}^{\Lambda,C}_T)$.
%Proposition~\ref{propcharacterizationbarXDexp}, and maximal %estimates
%given in Proposition~\ref{propLlogLinequality}
%Then, the running maximum of the entropic process $\rho_{\cdot}(
%%exponential of its norm is dominated by the Shannon entropy of $\exp(
%$\mathbb E(\exp(\bar{X}^{\Lambda,C}_T) \bar{X}^{\Lambda,C}_T)-
%e^{C_s}\,d\Lambda_s$ belongs to $\mathbb L^1_{\exp}$ and the
%%exponential of its norm is also dominated by the Shannon entropy of $
%defined as $N^0_t \equiv
%E(\exp{\bar{X}^{\Lambda,C}_T(|Y_T|)}|\cF_{\cdot} )$ is in
%$\mathbb{H}^1$. Moreover, both processes $\Phi$ and its running
%supremum $\max\Phi$ are also in $\mathbb{H}^1$.
%%%%%%%%
%================================
%s4 #&#
\section{Quadratic variation estimates and stability results}\label
{sectionestimates}
%================================
%================================
We are now capable to establish the main contribution of this paper,
that is, some
stability results, which require some uniform estimation of key
quantities, including quadratic variation and running supremum. In
order to use the previous inequalities, we need the family of
$\mathcal{Q}(\Lambda,C)$-semimartingales we consider to be uniformly
dominated. Following Remark~\ref{remetaT}, we can
replaced $Y_T$ by a generic random variable $\eta_T$ such that $|\eta
_T| \geq|Y_T|$ and $\bar X^{C,\Lambda}_T(|\eta_T|)$ satisfies an
appropriate integrability condition. Therefore, it seems natural to
introduce the following
class ${\mathcal S_Q}(|\eta_T|,\Lambda, C)$, and to work within this
class of quadratic semimartingales:
%de4.1 #&#
%
\begin{edefinition}\label{defSQ}
Let $|\eta_T|$ be a $\cF_T$-random variable, such that\break
$\bar X^{C,\Lambda}_T(|\eta_T|)=e^{C_{T}}|\eta_T|+\int
_0^Te^{C_{s}}\,d\Lambda_s$ belongs to ${\mathbb
L}^1_{\exp}$. The class ${\mathcal S_Q}(|\eta_T|,\Lambda, C)$ is the
set of $\mathcal{Q}(\Lambda,C)$-semimartingales $Y_{\cdot}$ defined on
$[0,T]$, such that
$|Y_{\cdot}|\leq
\rho_{\cdot}(e^{C_{\cdot,T}}|\eta_T|+\int_{\cdot}^Te^{C_{\cdot,s}}\,
d\Lambda
_s)$ a.s.
\end{edefinition}

%================================
%s4.1 #&#
\subsection{Quadratic variation estimates}\label{subvarqestimates}
%================================
We now study the quadratic variation of
$\mathcal{Q}(\Lambda,C)$-semimartingale $Y_{\cdot}$ when $Y_{\cdot}
$ belongs to
$\mathcal{S_Q}(|\eta_T|,\Lambda, C)$. Following Kobylanski
\cite{Kobylanski00}, the best way to do so
is to use the function $v(x)=e^{x}-1-x$ instead of the simple
exponential function. This function is indeed positive, convex, and increasing
for $x\geq0$, and satisfies $v''(x)-v'(x)=1$.
%Observe also that $v(x)=u(e^x)$, where $u$ is the function introduced
%in Proposition~\ref{propL%logLinequality} {\cgblue for $m=1$}. {
%consider the logarithm of exponential martingales}.
In the following, we use the short notation $\bar X_T^{C,\Lambda
}(|\eta
_T|)=\bar X^{C,\Lambda}_T$.
%
%th4.2 #&#
%
\begin{theorem}[(Quadratic variation estimates)]\label{thevarqestimates}
Let $Y_{\cdot} \in{\mathcal S_Q}(|\eta_T|,\Lambda, C)$.
\begin{longlist}[(iii)]
\item[(i)]Then, the quadratic variation $\varqM_{\cdot}$ of the $\cQ
(\Lambda,
C)$-semimartingale $Y_{\cdot}=Y_0+M_{\cdot}-V_{\cdot}$ satisfies for
any stopping times
$\sigma\leq T$,
%e19 #&#
%
\begin{equation}\label{eqquadraticvariationestimate}
\quad\tfrac{1}{2}\E[\varqM_{\sigma,T}|\cF_\sigma]\leq\Phi_\sigma
(|Y_T|){\mathbf
1}_{\{\sigma<T\}}\leq\E\bigl[\exp(\bar X^{C,\Lambda}_T(|\eta_T|))
{\mathbf1}_{\{\sigma<T\}}|\cF_\sigma\bigr].
\end{equation}
In particular, the martingale $M_{\cdot}$ is in ${\mathbb H}^{2}$,
\textit{with
the uniform estimate}
%a uniform control of the quadratic norm
%e20 #&#
%
\begin{equation}\label{eqquadraticvariationH2estimate}
\E\bigl[\tfrac{1}{2}\varqM_T\bigr]\leq\E[\exp(\bar X^{C,\Lambda}_T(|\eta_T|))].
\end{equation}
\item[(ii)]Let $p^\eta=\sup\{p; \E[\exp(p \bar
X^{C,\Lambda}_T(|\eta_T|))]<+\infty\}$. Then $p^\eta\geq1$ and
$\forall p\in[1,
p^\eta[$, the martingale $M$ belongs to ${\mathbb H}^{2p}$, and
%e21 #&#
%
\begin{equation}\label{eqquadraticvariationHpestimate}
\E[\varqM_T^p]\leq(2p)^p \E[\exp(p \bar X^{C,\Lambda}_T(|\eta
_T|))].
\end{equation}
\item[(iii)]If $\Phi_t(|\eta_T|)=\E[\exp
(e^{C_{t,T}}|\eta_T|+\int_t^Te^{C_{t,u}}\,d\Lambda_u)|\cF_t]$ is
uniformly bounded in $t \leq T$, then the conditional quadratic
variation $\frac{1}{2}\E[\varqM_{\sigma,T}|\cF_\sigma]$ is uniformly
bounded. Hence $M_{\cdot}$ is a BMO-martingale.
\end{longlist}
\end{theorem}

\begin{pf}
%In the proof, for the sake of simplicity in the notation, the
%dependency of $X^C_T$ and $\Phi$ on respect to $|\eta_T|$ is not
%mentioned.\\
By analogy with the previous notation, when using the
function $v(x)=e^{x}-1-x$, we set $V^{\Lambda,C}_t(e^{|Y|})=
v(|Y_t|)+\int_0^t v'(|Y_s|)(d\Lambda_s+|Y_s|\,dC_s)= v(|Y_t|)+\int_0^t
v'(|Y_s|)\,dD^{\Lambda,C}_s$. So,\vspace*{1pt} $U^{\Lambda
,C}_t(e^{|Y|})-V^{\Lambda
,C}_t(e^{|Y|})=1+|Y_t|+D^{\Lambda,C}_t(|Y|)$, and both processes
$U^{\Lambda,C}_{\cdot}$ and $V^{\Lambda,C}_{\cdot}$ are in the
class $(\cD)$ since
$Y_{\cdot} \!\in\!{\mathcal S_Q}(|\eta_T|,\Lambda, C)$.

\begin{longlist}[(iii)]
\item[(i.1)] As we see in the proof of Proposition \ref
{propcharacterizationbarXDexp},
the semimartingale $|Y_{\cdot}|$ is associated
with the martingale $M^s_{\cdot}=\operatorname{sign}(Y_{\cdot})\star M_{\cdot
}$, the finite variation
process $V^s_{\cdot}=\operatorname{sign}(Y_{\cdot})\star V_{\cdot}$ and the
local time at $\{0\}$,
that disappears in the It\^o's formula since $v'(0)=0$. Using similar
calculation to those of the previous section, and the identity
$v''(x)-1=v'(x)$, we obtain that the process
$V^{\Lambda,C}_t(e^{|Y|})-\frac{1}{2}\varqM_t =v(|Y_0|)+\int_0^t
v'(|Y_s|)\,dM^s_s+\int_0^t v'(|Y_s|)(dD^{\Lambda,C}_s-dV^s_s+\frac{1}{2}
\,d\varqM_s)$
is a submartingale, and since $V^{\Lambda,C}_{\cdot}$ is in the class
$(\cD)$,
for any $\sigma\leq T$,
$ \E[\frac{1}{2}\varqM_{\sigma,T}|\cF_\sigma]\leq\E
[v(|Y_T|)-v(|Y_\sigma
|)+\int_\sigma^Tv'(|Y_s|)\,dD^{\Lambda,C}_s|\cF_\sigma]$.

\item[(i.2)] Since, by definition, $\forall x\geq0, 0\leq v(x) \leq e^x$ and $
v'(x) \leq e^x$,
\[
\int_\sigma^Tv'(|Y_s|)\,dD^{\Lambda,C}_s \leq\int
_\sigma^T
\Phi_s(d\Lambda_s+\ln|\Phi_s|\,dC_s) \qquad\mbox{for any }\sigma\leq T.
\]

Thanks to the supermartingale property of $U^{\Lambda,C}_{\cdot}(\Phi)$
[Proposition~\ref{propUtransform}(ii)] and the inequality $\Phi
_{\cdot}\geq
\exp(|Y_{\cdot}|)$ [implying in particular $\Phi_T=\exp(|\eta
|)\geq
v(|Y_T|)$], we have
$ \E[\int_\sigma^T \Phi_s(d\Lambda_s+\ln|\Phi_s|\,dC_s)|\cF
_\sigma]\leq\E
[\Phi_\sigma-\Phi_T|\cF_\sigma]$ and
\begin{eqnarray*}
\E\bigl[\tfrac{1}{2}\varqM_{\sigma,T}|\cF_\sigma\bigr]&\leq&\E
[v(|Y_T|)-v(|Y_\sigma
|)-(\Phi_T-\Phi_\sigma)|\cF_S]\\
&=&\E\bigl[\bigl(-\bigl(\Phi_T-v(|Y_T|)+v(|Y_\sigma|)\bigr)+\Phi_\sigma\bigr)\mathbf
1_{\{\sigma<T\}}\big|\cF_\sigma\bigr]\\
&\leq&\Phi_\sigma\mathbf1_{\{\sigma<T\}}\leq\E\bigl[\exp\bar
X^{\Lambda
,C}_{T} {\mathbf1} _{\{\sigma<T\}}\big|\cF_\sigma\bigr].
\end{eqnarray*}
%
%
% \\[2mm]
%
\item[(ii)]As observed in Lenglart, L\'epingle and Pratelli \cite
{Lenglart-Lepingle-Pratelli}, the final result is a simple consequence
of the
so-called Garsia--Neveu lemma (Lemma~\ref{lemNeveuGarsia}) (see, e.g., Neveu~\cite{Neveu}) recalled below.

\item[(iii)]This is a straightforward consequence of the inequality
$\E[\frac{1}{2}\varqM_{\sigma,T}|\break\cF_\sigma] \leq\Phi_\sigma
(|\eta_T|)$.\quad\qed
\end{longlist}
\noqed\end{pf}
%
%le4.3 #&#
%
\begin{elemma}[(Garsia--Neveu lemma)]
\label{lemNeveuGarsia} Let $A_{\cdot}$ be a predictable c\`adl\`ag
increasing process and $U$ a random variable, positive and integrable.
If for any stopping times $\sigma\leq T$, $\E[A_T-A_\sigma|\cF
_{\sigma}]\leq\E[U{\mathbf1} _{\{\sigma<T\}}|\cF_{\sigma}],$
\[
\forall p\geq1\qquad \E[A_T^p]\leq p^p\E[U^p].
\]

More generally, $\E[F(A_T)]\leq\E[F(pU)]$ for any convex function $F$
such that $p=\sup_{x>0}(x (\ln F)'(x))<+\infty$.
\end{elemma}

% \centerline{$\E[A_T-A_\sigma|\cF_{\sigma}]\leq\E[U{\mathbf1} _{\{
% {\mathbf1} _{\{0<\sigma\}}
% \end{remark}

Here we apply this lemma to the random variable $U=\exp( \bar
X^{\Lambda,C}_{T}(|\eta_T|))$
for any $p\geq1$ such that $U\in{\mathbb L}^p.$
As a corollary of this result, uniform estimates may be obtained for
the total variation of the process~$V_{\cdot}$.
%
%co4.4 #&#
%
\begin{ecorollary}\label{corVestimates}
Let $Y_{\cdot}\in{\mathcal S_Q}(|\eta_T|,\Lambda, C)$.
The total variation of the process~$V_{\cdot}$ such that $Y_{\cdot
}=Y_0+M_{\cdot}-V_{\cdot}$
satisfies for $1\leq p<p^\eta$
%e22 #&#
%
\begin{equation}\label{eqquadraticvariationHpestimate}
\E[|V |^p_T]\leq(2p)
^p \E[\exp(p \bar X^{C,\Lambda}_T)].
\end{equation}
When $\Phi_{\cdot}(|\eta_T|)=\E[\exp
(e^{C_{\cdot,T}}|\eta_T|+\int_{\cdot}^Te^{C_{t,u}}\,d\Lambda_u)|\cF
_{\cdot}]$ is
bounded by $K_C$, then $\mathbb{E}[|V|_{\sigma
,T}|\cF_\sigma]\leq2K_C.$
%{\cgblue revoir les hypotheses minimales - voir ce que l'on fait a
%ce sujet}
\end{ecorollary}

\begin{pf} Since $V_{\cdot}$ satisfies the structure condition $\cQ
(\Lambda,
C)$, $\E[ |V|_{\sigma,T}|\cF_\sigma]\leq
E[\Lambda_{\sigma,T}+\int_\sigma^T|Y_s|\,dC_s+\frac{1}{2}\varqM
_{\sigma
,T}|\cF
_\sigma]\leq
2\E[\exp(\bar X^{\Lambda,C}_{T}){\mathbf1} _{\{\sigma<T\}}|\cF
_\sigma]$.
Indeed,\break $\E[\Lambda_{\sigma,T}+\int_\sigma^T|Y_s|\,dC_s|\cF
_\sigma
]\leq\E
[\int_\sigma^T
e^{|Y_s|}(d\Lambda_s+|Y_s|\,dC_s)|\cF_\sigma] \leq
\E[(\Phi_\sigma-\Phi_T)| \cF_\sigma] \leq\E[\exp( \bar
X^{\Lambda,C}_{T}){\mathbf1} _{\{\sigma<T\}}|\cF_\sigma]$. We
conclude with Lemma~\ref{lemNeveuGarsia}.
\end{pf}
%
%====================================
%s4.2 #&#
\subsection{\texorpdfstring{Stability results for $\mathcal{Q}(\Lambda,C)$-semimartingales}
{Stability results for Q(Lambda,C)-semimartingales}}\label{substabilityresults}\label{substabilityresults}

%====================================
We can start by noticing that the class ${\mathcal S_Q}(|\eta
_T|,\Lambda
, C)$ is stable by a.s. convergence,
since the submartingale property of both processes $U_{\cdot}(e^Y)$ and
$U_{\cdot}(e^{-Y})$, dominated by the $(\cD)$-supermartingale
$U_{\cdot}(\Phi)$, is
stable by a.s. convergence.
Moreover, Theorem~\ref{thcharacterizationgeneralquasimart} implies that\vadjust{\goodbreak}
the limit process is continuous and is also in ${\mathcal S_Q}(|\eta
_T|,\Lambda, C)$.
However, previous estimates of both
quadratic variation and
finite variation processes suggest that a better stability result
may hold true, in particular regarding the strong convergence of the
martingale parts. The space of martingales, where this convergence
takes place, depends essentially on the exponential integrability
properties of the random variable $X^{\Lambda,C}_T(|\eta_T|)$. The
method is very similar to that of Lepeltier and San Martin~\cite
{LepeltierSanMartin97}. When the
$\cQ(\Lambda,C)$-semimartingales are bounded, this type of results
has already been obtained for the $\mathbb H^2$-convergence by
Kobylanski~\cite{Kobylanski00} and Morlais~\cite{Morlais1}.
Our stability result is novel and direct, and
gives better convergence results with the ${\mathbb H}^1$ convergence.
This result, that appears here for the first time in a BSDE framework,
is based on an old result of Barlow and Protter~\cite{Barlow-Protter} on
the convergence of semimartingales.
%
%th4.5 #&#
%
\begin{theorem}\label{thNicolas}
Assume the sequence $(Y^n_{\cdot})$ of ${\mathcal S_Q}(|\eta
_T|,\Lambda, C)$
semimartingales is a Cauchy sequence for the
a.s. uniform convergence, that is,\break $\sup_{t\leq T}|Y^n_t-Y^{n+p}_t|$
tends to 0 almost surely when $n \rightarrow\infty$.
Then the limit process $Y_{\cdot}$ is a ${\mathcal S_Q}(|\eta
_T|,\Lambda,
C)$-semimartingale
$Y_{\cdot}=Y_0+M_{\cdot}-V_{\cdot}$.

Different types of convergence hold true for the processes $(M^n_{\cdot},
V^n_{\cdot})$ of the decomposition $Y^n_{\cdot}=Y^n_0+M^n_{\cdot
}-V^n_{\cdot}$:
\begin{enumerate}[(ii)]
\item[(i)]Martingales convergence of $(M^n_{\cdot})$ to $M_{\cdot}$.

\begin{enumerate}[(a)]
\item[(a)] The sequence $(M^n_{\cdot})$ converges to~$M_{\cdot}$ in
${\mathbb H}^1$.

\item[(b)]If for some $p>1$ $\bar X^{\Lambda,C}_{T}(|\eta_T|)\in
{\mathbb
L}^p_{\exp}$, the sequence $(M^n_{\cdot})$ converges to~$ M_{\cdot}$ in $\mathbb H^{2p}$, and in the BMO-space if $ \Phi
_S(|\eta
_T|)$ is bounded.
\end{enumerate}
%
%$(M^n_{\cdot})$ converges to a martingale $M_{\cdot}$ in the BMO-space.
\item[(ii)]The sequence of finite variation processes $(V^n_{\cdot})$
converges at
least in $\mathcal S^1$ to the process $V_{\cdot}$ satisfying the structure
condition $\cQ(\Lambda,C)$.
\end{enumerate}
\end{theorem}

%%%%%
%
\begin{pf} We proceed\footnote{An earlier proof of this result in the
BMO case is due to Nicolas Cazanave, a former Ph.D. student at Ecole
Polytechnique.} in several steps to prove this convergence result.
We first introduce some notation and make some elementary
calculations. For $s\leq t$, let $Y^{i,j}_t=Y_t^i-Y^j_t,
M^{i,j}_t=M_t^i-M^j_t $ and
$Y^{i,j}_{s,t}=(Y_t^i-Y^i_s)-(Y_t^j-Y^j_s) $, and the short notation
first introduced in Section~\ref{parmaximalinequalities}, $\sup
_{s\leq u\leq t}
| Y_{u}^{i,j}-Y_{s}^{i,j}|\equiv\max| Y^{i,j}_{s,t}|$.
Then for any stopping times $\sigma\leq\tau\leq T$,
\begin{eqnarray*}
\varq{M^{i,j}}_{\sigma,\tau}&=&|Y^{i,j}_{\sigma,\tau}|^2-2\int
_\sigma
^\tau Y^{i,j}_{\sigma,s} \,dY^{i,j}_s \\
&\leq&|Y^{i,j}_{\sigma,\tau}|^2-2\int_\sigma^\tau Y^{i,j}_{\sigma
,s}
\,dM^{i,j}_s
+ 2\int_\sigma^\tau|Y^{i,j}_{\sigma,s}| \, d(|V^j|_s+|V^i|_s).
\end{eqnarray*}
Using either the fact that $Y^{i,j}$ is bounded, or a uniform
localization procedure, the stochastic integral
$\int_\sigma^{\tau_n} Y^{i,j}_{\sigma,s} \,dM^{i,j}_s$ has null
conditional expectation for\vadjust{\goodbreak} a well-chosen stopping time $\tau_n$. Then,
thanks to the monotonicity of $\varq{M}$ and Corollary~\ref
{corVestimates}, with $B^{i,j}=2(|V^i|+|V^j|)$,
\begin{eqnarray*}
\mathbb{E}[\varq{M^{i,j}}_{\sigma,T} | \cF_\sigma]
&\leq&\mathbb{E}\biggl[\max|Y^{i,j}_{\sigma,T}|^2 {\mathbf1} _{\{\sigma
<T\}}+
\int_\sigma^T\max|Y^{i,j}_{\sigma,s}| \,dB^{i,j}_s \Big| \cF_{\sigma}
\biggr]\\
&\leq&\mathbb{E}\bigl[(\max|Y^{i,j}_{0,T}|^2+\max
|Y^{i,j}_{0,T}|B^{i,j}_T) {\mathbf1} _{\{\sigma<T\}}\big| \cF_{\sigma
}\bigr].
\end{eqnarray*}
We now start with the proof corresponding to the assumption $\bar
X^{\Lambda,C}_{T}(|\eta_T|)\in{\mathbb L}^p_{\exp}$, since it is very
similar to the proof in the linear growth case (see Lepeltier and San
Martin~\cite{LepeltierSanMartin97}).

(i.b) Thanks to the Garsia--Neveu lemma (Lemma \ref
{lemNeveuGarsia}), for $r\geq1$,
\begin{eqnarray*}
\mathbb{E}[\varq{M^{i,j}}_T^r ] &\leq& r^r \mathbb{E}[(\max
|Y^{i,j}_{0,T}|^2+\max|Y^{i,j}_{0,T}| B^{i,j}_T)^r]\\
& \leq&\tfrac{1}{2}(2r)^r\{\mathbb{E}[(\max|Y^{i,j}_{0,T}|)^{2r}
]+\mathbb{E}[(\max|Y^{i,j}_{0,T}| B^{i,j}_T)^r]\}.
\end{eqnarray*}
Then, since $B^{i,j}_T$ belongs to $\mathbb L^p$, by H\"older
inequalities, for any $p$ and $q$ such that $\frac{1}{p}+\frac
{1}{q}=1$, and $1\leq r<p$, if $K_r=\frac{1}{2}(2r)^r $
\begin{eqnarray*}
\mathbb{E}[(\max|Y^{i,j}_{0,T}| B^{i,j}_T)^r]&\leq&(\mathbb
{E}[(\max|Y^{i,j}_{0,T}|)^q])^{{r}/{q}}
(\mathbb{E}[(B^{i,j}_T)^p])^{{r}/{p}},\\
\mathbb{E}[\varq{M^{i,j}}_T^r ]&\leq&
K_r\{\mathbb{E}[\max|Y^{i,j}_{0,T}|^{2r}
]\\
&&\hspace*{15pt}{}+(\mathbb{E}[(\max|Y^{i,j}_{0,T}|)^q])^{{r}/{q}}
(\mathbb{E}[(B^{i,j}_T)^p])^{{r}/{p}}\}.
\end{eqnarray*}
From the monotonicity of both sides of this inequality with respect to
$r$, we can take $r=p$.
%{\cgreen Faire le lien avec ce qu'on a raconte sur le llogl}
We have used that
$\max|Y^{i,j}_{0,T}|$ has finite moments of all orders since as shown
in Section~\ref{parmaximalinequalities} $\max|Y^i_{0,T}|$ and $\max
|Y^i_{0,T}|$ are in ${\mathbb L}^p_{\exp}$.
Hence, we have the desired convergence.

(i.c) In the bounded case, thanks to Corollary~\ref{corVestimates},
the conditional total variation $\mathbb{E}[|V^n|_{\sigma,T}|\cF
_\sigma
]$ are
uniformly bounded by $C_V$. To obtain the BMO convergence, we have to modify
the previous proof, by using an integration by parts formula
involving the conditional variation of $B^{i,j}$,
\begin{eqnarray*}
\mathbb{E}\biggl[\int_\sigma^T\max| Y_{\sigma,s}^{i,j}| \,dB^{i,j}_s |
\mathcal{F}_\sigma\biggr]
&=& \mathbb{E}\biggl[\int_\sigma^T d_u\max| Y_{\sigma,u}^{i,j}|\biggl(
\mathbb{E}\biggl[\int_u^T \,dB^{i,j}_s \Big| \mathcal{F}_u\biggr]\biggr) \Big| \mathcal
{F}_\sigma\biggr]\\
&\leq&
2 C_V
\mathbb{E}[\max| Y_{\sigma,T}^{i,j}| | \mathcal{F}_\sigma],
\end{eqnarray*}
and so
\[
\mathbb{E}[\varq{M^{i,j}}_{\sigma,T} |
\mathcal{F}_\sigma]\leq2 C_V
\mathbb{E}[\max| Y_{\sigma,T}^{i,j}| | \mathcal{F}_\sigma
]+\mathbb
{E}[|Y^{i,j}_{\sigma,T}|^2 | \cF_\sigma].
\]

%In terms of quadratic variation, we have:
%$\mathbb{E}[\varq{M^{i,j}_{\sigma,T}} | \mathcal{F}_\sigma]\leq2
% C_V
%[|Y^{i,j}_{t,T}|%^2 | \cF_t].$
Then, the BMO-convergence holds true.

(i.a) The proof of the general case requires a different argument,
based on a result of Barlow and Protter~\cite{Barlow-Protter}
on the convergence of semimartingales. In the framework of quadratic
semimartingales, the key points are the uniform estimates of both
the quadratic variation and the total variation\vadjust{\goodbreak} given in Theorem \ref
{thevarqestimates}, equation \eqref{eqquadraticvariationH2estimate}
and Corollary~\ref{corVestimates}. The proof given in \cite
{Barlow-Protter} of the $\mathbb H^1$-convergence of the martingales is
based on the
square root of the inequality given at the beginning of the proof,
\[
\varq{M^{i,j}}_{t}\leq|Y^{i,j}_{0,t}|^2-2\int_0^t Y^{i,j}_{0,s} \,dM^{i,j}_s+
2\int_0^t |Y^{i,j}_{0,s}| \,dB^{i,j}_s.
\]

The first step is to estimate the square root of $\max
|Y^{i,j}_{0,\cdot}\star M^{i,j}_{0,T}|$
%{\cgreen Pauline: j'ai des questions de notation ici et par la suite.
%Doit-on utiliser $\star$ ou $\star$? Pour les IS on utilise \star,
%pour les proc variation finies l'asterisque}
using the Burkolder--Davis--Gundy inequalities
for continuous martingales for $p=\frac{1}{2}$, that have been
recalled in
Section~\ref{subsecentropicinequalities}:
$\E[\max|Y^{i,j}_{0,\cdot}\star M^{i,j}_{T}|^{{1}/{2}}]\leq
\bar{C}
\E[\varq{Y^{i,j}_{0,\cdot}\star M^{i,j}}_{T}^{1/4}]$ where $\bar
{C}$ is a
universal constant. Then, since $\E[\varq{Y^{i,j}_{0,\cdot}\star
M^{i,j}}_{T}^{1/4}]\leq\E[(\max|Y^{i,j}_{0,T}|)^{{1}/{2}}\varq
{M^{i,j}}_{T}^{1/4}]$,
\begin{eqnarray*}
\E\bigl[\sqrt{\varq{M^{i,j}}_T}\bigr]& \leq&\E[\max|Y^{i,j}_{0,T}|]+\sqrt
{2}\bar{C} \E[\max|Y^{i,j}_{0,T}|]^{ {1}/{2}}\E\bigl[\sqrt{\varq
{M^{i,j}}_T}\bigr]^{{1}/{2}}\\
&&{}+
\sqrt{2}\E[\max|Y^{i,j}_{0,T}|]^{{1}/{2}}\E[B^{i,j}_T|]^{{1}/{2}}.
\end{eqnarray*}
Since $\E[\sqrt{\varq{M^{i,j}}_T}] $ and $\E[B^{i,j}_T]$ are uniformly
bounded, and $\E[\max|Y^{i,j}_{T}|]$ goes to~$0$, then $\E[\sqrt
{\varq{M^{i,j}}_T}] $ also goes to $0$. The $\mathbb
H^1$-convergence of the martingale part is established.

(ii) The next point is to study the convergence of the sequence
$(V^n_{\cdot})$ to a process $V_{\cdot}$ satisfying the same
structure condition
$\mathcal{Q}(\Lambda,C)$. Since, the sequence $(Y^n_{\cdot},
M^n_{\cdot},
\varq{M^n}_{\cdot}^{{1}/{2}})$ converges in $\mathcal S^1$ to
$(Y_{\cdot}, M_{\cdot},\varqM^{{1}/{2}}_{\cdot})$, the sequence
$(V^n_{\cdot})$ also converges in
$\mathcal S^1$. Therefore, we can extract a
subsequence, still denoted $(Y^n_{\cdot}, M^n_{\cdot},V^n_{\cdot},$
$\varq{M^n}_{\cdot}^{{1}/{2}})$, such that the
sequence converges uniformly in time almost surely.

(iii) This point is obvious since, as observed at the beginning of
this section, the class ${\mathcal S}(|\eta_T|)$ is stable by a.s. convergence.
\end{pf}
%
%Thanks to the \textit{structure condition}, the uniformly dominated
%increasing processes
%$\overline{A}^n_{\cdot}=-V^n_{\cdot}+\demi\varq{M^n}_{\cdot}+\Lambda_{
%C_{\cdot}$
%and
%$\underline{A}^n_{\cdot}=V^n_{\cdot}+\demi\varq{M^n}_{\cdot}+\Lambda_{
%C_{\cdot}$
%also converge uniformly almost surely to the increasing processes
%$\overline{A}_{\cdot}$ and $\underline{A}_{\cdot}$. The process
%$V_{\cdot}=\demi(\overline{A}_{\cdot}-\underline{A}_{\cdot})$, which
%is the difference
%of two increasing processes summing up to
%$\demi\varq{M}_{\cdot}+\Lambda_{\cdot}+|Y|_{\cdot}\star C_{\cdot}$,
%also satisfies
%the structure condition. \\
%Note that we could use the characterization of
%$\mathcal{Q}(\Lambda,C)$-semimartingales given in Theorem
%property.
%{We have made the choice of an unified presentation of convergence
%results for uniformly convergente sequence of $\mathcal{Q}(
%assumption \rmi b) on the structure of $\bar X^{\Lambda,C}_{T}(|
%convergence of the sequence.}
% \begin{ecorollary} Assume $\bar X^{\Lambda,C}_{T}(|\eta_T|)\in{
%belongs to $
%%%%%%%%%%%%%%%%%%%%%%%%%%%%%%%%%
%pa4.2.0.1 #&#
\subsubsection*{Stability results for BSDE-like quadratic
semimartingales}
%%%%%%%%%%%%%%%%%%%%%%%%%%%%%%%%%%%%%%
%{\cgreen Pauline: Je trouve cette section tres claire maintenant.}\\
To obtain the convergence of the finite variation processes in total
variation, we need to make additional assumption on the processes
$V^n$, as in the BSDE framework. We adopt the general setting where the
reference to the Brownian framework is relaxed as in El Karoui and Huang
\cite{NEKHuang}.
%precise the form of the processes $V^n$, as in the BSDE framework,
%where the %reference to the Brownian setting is relaxed as in El
%Karoui and Huang \cite%{NEKHuang}.
%de4.6 #&#
%
\begin{edefinition}[(BSDE-like quadratic semimartingale)]\label
{defgeneralquadraticBSDE}
Let us consider a continuous predictable increasing process $K_{\cdot}$,
%a general BSDE framework with
a $d$-dimensional continuous orthogonal
martingale $N_{\cdot}=(N^i_{\cdot})_{i=1}^d$, with quadratic
variation $\varq
{N^i}_{\cdot}$ strongly dominated by $K_{\cdot}$ such that $d\varq
{N^i}_t=\gamma
^{i}_t \,dK_t$,
%the diagonal predictable quadratic
%variation matrix is dominated by some continuous predictable
%increasing process $K_{\cdot}$, with $d\varq{N^i}_t=\gamma^{i}_t dK_t$.
two increasing processes $\Lambda_{\cdot}$ and $C_{\cdot}$ also
dominated by $K_{\cdot}$
such that $d\Lambda_t=l_t\,dK_t$ and $dC_t=c_t \,dK_t$, such that all processes
$\gamma^{i}_{\cdot}, l_{\cdot},c_{\cdot}$ are bounded by $k$ (e.g., $K_{\cdot}=\sum
_{i=1}^d\varq{N^i}_{\cdot}+\Lambda_{\cdot}+C_{\cdot}$ and $k=1$).
The coefficient
$g(\cdot,y,z)$ is a $\mathcal{P}\otimes\mathcal{B}(\mathbb{R}\times
\mathbb
{R}^d)$ measurable process, often assumed to be continuous with respect
to $(y,z)$.\vadjust{\goodbreak}

A semimartingale $Y_{\cdot}$ with the decomposition $Y_{\cdot
}=Y_0-V_{\cdot}+M_{\cdot}$ is said
to have a quadratic coefficient $g$ if $dY_t=-dV_t+dM_t,$ with
%e23 #&#
%
\begin{equation}\label{eqgeneralquadraticBSDE}
\cases{
% \begin{array}{cllll}
%
dV_t=g(t,Y_t, Z_t) \,dK_t,\vspace*{2pt}\cr
dM_t=Z_t\,dN_t+dM_t^\bot \qquad\forall i\,
d\varq{N^i,M^\bot}_t=0, \vspace*{2pt}\cr
\displaystyle|g(t,Y_t, Z_t)| \leq\frac{1}{\delta}l_t +|Y_t| c_t +\frac{\delta
}{2} \bigl|\sqrt\gamma_t Z_t\bigr|^2,\vspace*{2pt}\cr
\displaystyle\bigl|\sqrt\gamma_t Z_t\bigr|^2=\sum
_{i=1}^d\gamma^i_t |Z^i_t|^2.
}
\end{equation}
The local martingale $Z \star N$ is the orthogonal projection of the local
martingale~$M_{\cdot}$ onto the space of stochastic integrals
generated by
the local martingale $N_{\cdot}$, and $d\varq{Z \star N}_t\ll \,d\varqM
_t$, so
that $d|V|_t\ll\frac{1}{\delta}\,d\Lambda_t+|Y_t| \,dC_t +\delta
\,d\varqM
_t$ and $Y_{\cdot}$
is a quadratic semimartingale.
\end{edefinition}

When considering sequences of BSDE-like quadratic semimartingales
under mild assumptions on the sequence of coefficients, the sequence
of finite
variation processes is converging in total variation in the
appropriate space, and the limit is still a BSDE-like quadratic
semimartingale.

The uniform convergence of the quadratic semimartingales
needed for these convergence results may seem very strong. We know
however from Theorem~\ref{thcharacterizationgeneralquasimart}
that all the processes obtained by a.s. convergence are continuous.
Thanks to Dini's theorem, the monotone convergence implies uniform
convergence for continuous functions on compact spaces.
Therefore, by a localization procedure, we can prove
the following very strong result.
%
%th4.7 #&#
%
\begin{theorem}\label{thmonotonecv}
Let assume the sequence $(Y^n_{\cdot})$ to be a monotone sequence of
${\mathcal S_Q}(|\eta_T|,\Lambda, C)$-semimartingales
%quadratic
converging almost surely to a process $Y_{\cdot}$.
\begin{enumerate}[(ii)]
\item[(i)]Then, the limit process $Y_{\cdot}$ is a continuous ${\mathcal
S_Q}(|\eta
_T|,\Lambda, C)$-
%quadratic
semimartingale,
the convergence is locally uniform and all properties given in Theorem
\ref{thNicolas} hold (locally) true. In particular, there exists a
subsequence of martingales $M^n_{\cdot}=Z^n\star N_{\cdot}+M^{n,\bot
}_{\cdot}$ converging
in $\mathbb{H}^1$ and almost surely to $M_{\cdot}=Z\star N_{\cdot
}+M^{\bot}_{\cdot}$.

\item[(ii)]Suppose in addition that the processes $(Y^n_{\cdot})$ are BSDE-like
quadratic semimartingales, associated with a sequence of monotone
coefficients $g_n$ converging almost surely to $g$, having the
following properties:
%such that $Y^n_{\cdot}=Y^n_0+Z^n\star N_{\cdot}+M^{n,\bot}-g_n(\cdot,Y^n_{
%K_{\cdot}$
%and for which the following assumptions are made:
\begin{enumerate}[(a)]
\item[(a)]The monotone sequence $g_n$ have uniform quadratic growth:
\[
|g_n(t,Y^n_t,Z^n_t)| \leq\frac{1}{\delta}l_t +|Y^n_t|
c_t +\frac{\delta}{2} \bigl|\sqrt\gamma_t Z^n_t\bigr|^2,\qquad  d\mathbb P\times
dK_{\cdot} \mbox{ a.s.}
\]

\item[(b)]The sequence $g_n(\cdot,Y^n_{\cdot},Z^n_{\cdot})$ converges to
$ g(\cdot,Y_{\cdot},Z_{\cdot}), d\mathbb P\times dK_{\cdot} \mbox{ a.s}$.
\end{enumerate}
Then, the limit process $Y_{\cdot}$ is a BSDE-like ${\mathcal
S}_Q(|\eta
_T|,\Lambda, C)$-semimartingale with coefficient $g(t,y,z)=\lim g_n(t,y,z)$.
\end{enumerate}
\end{theorem}

\begin{pf}
Note the characterization of ${\mathcal S}_Q(|\eta_T|,\Lambda,
C)$-semimartingales given in Theorem \ref
{thcharacterizationgeneralquasimart} passes to the limit, since all
processes $U^{\Lambda, C}_{\cdot}(e^{|Y^n|})$ are dominated by the
$(\cD)$-process $U^{\Lambda, C}_{\cdot}(\Phi(|\eta_T|))$. The
limit process
$Y$ is a continuous ${\mathcal S}_Q(|\eta_T|,\Lambda, C)$-semimartingale,
with decomposition $Y_{\cdot}=Y_0+M_{\cdot}-V_{\cdot}.$
\begin{longlist}[(ii)]
\item[(i)] The localization procedure is based on the family $(T_K)$ of
stopping times as to bound the u.i. martingale $N^0_t=\mathbb
E[\exp(\phi_0(|\eta_T|)|\cF_t]$ by $K$. By the characterization of
u.i. continuous martingale (see, e.g., Azema, Gundy and Yor
\cite{Azema-Gundy-Yor}), the sequence $T_K$ goes to $\infty$ and for
$K\geq K_\varepsilon$ large enough, $\mathbb P(T_{K}<T)\leq
\frac{\varepsilon}{K}$. Therefore, the sequence $(Y^n_{\cdot\wedge T_K})$
lives on a compact set where the monotone convergence to a
\textit{continuous process} is uniform.
The sequence of martingales $(M^n_{\cdot\wedge T_K})_n$ strongly
converges in the appropriate space to the martingale $M_{\cdot\wedge
T_K}$. The same property holds true for the sequence $V^n_{\cdot\wedge
T_K}$.
Thanks to the previous estimates, for all these processes $Y^n_{\cdot
},M^n_{\cdot},
V^n_{\cdot}$ the convergence is uniform on $[0,T\wedge T_K] $ in
probability.

\item[(ii)] Let $Z^{n,K}_t \equiv Z^n_t \ind_{\{t\leq T_K\}}$ in such way
that $(Z^{n}\star N)_{\cdot\wedge T_K}=Z^{n,K}\star N_{\cdot}$. Since the
sequence $(M^n_{\cdot\wedge T_K})_n$ strongly converges, the
sequences of orthogonal martingales $(M^{n,\bot}_{\cdot\wedge T_K})_n$
and $(Z^{n,K}\star N_{\cdot})_n$ also strongly converge in the
appropriate space, and at least in $\mathbb H^1$.

Therefore, we can extract a subsequence still denoted $Z_{\cdot}^{n,K}$
converging a.s. By assumption, for $t\leq T_K$ the sequence
$g^n(t,Y^n_t, Z^{n,K}_t)$ goes to $g(t,Y_t, Z_t)$ $dK_{\cdot}\otimes
d\mathbb
P$ a.s.
It now remains to show that the convergence is also true in expectation.
%$\mathbb E[\int_0^{T_K}|g_n(s,Y^{n}_s,Z^{n}_s-g(s,Y_s,Z_s)|\,dK_s]$
%goes
%to $0$.\\
Observe that $\mathbb
E[\int_0^{T_K}|g_n(s,Y^{n}_s,Z^{n})_s-g(s,Y_s,Z_s)|\ind_{\{|Z^{n}_s|
\leq C\}}\,dK_s]$ goes to $0$, by
dominated convergence, since $\Phi_{\cdot}$ and $Y_{\cdot}^n$ are
bounded on $[0,T_K]$.
Moreover, since the sequence in $n$ of
the quadratic variations at time $T_K$, $\varq{Z^{n,K}\star N}_{T_K}$
is bounded in~$\mathbb
L^1$, for $s\leq T_K$, $|g_n(s,Y^{n}_s, Z^{n}_s)-g(s,Y_s,Z_s)| \leq
\Psi
_s+\frac{1}{2}
|Z^n_s|^2$,
with $\Psi_t \ind_{\{t\leq T_K\}} \in\mathbb L^1(d\mathbb P\otimes
dK_s)$ and
$\mathbb P(|Z^{n}_s| \geq C)\leq\frac{1}{C^2}\mathbb E
(|Z^n_s|^2)$. Hence, $\mathbb
E[\int_0^{T_K}|g_n(s,Y^{n}_s, Z^{n}_s)-g(s,Y_s, Z_s)|\ind_{\{
|Z^{n}_s| >
C\}}\,dK_s]$ goes to $0$ when $C$ goes to $\infty$, uniformly in $n$.
As a consequence, the process $V_{\cdot}$ in the decomposition of the quadratic
semimartingale~$Y_{\cdot}$ is given by $dV_t=g(t,Y_t, Z_t) \,dK_t$ on $[0,T_K]$
for any $K$.\quad\qed
\end{longlist}
\noqed\end{pf}
%
%re4 #&#
%
\begin{remark} Delbaen, Hu and Bao show in~\cite{Delbaen-Hu-Bao} that
increasing the growth of the coefficient into a superquadratic growth
yields to ill-posed problems. In particular, monotone stability does
not hold any more. For classical BSDEs, when the coefficient simply
depends on $z$, superquadratic growth means that $\limsup
g(z)/|z|^2=\infty$.
\end{remark}

%%%%%
%%%%%%%%%%%%%%%%%%%%%%%%
%Most of the papers in the literature focusing on the study
%of quadratic BSDEs consider the situation where the martingale $M_{
%BMO, as this gives a well-known framework for the existence of a
%solution for the BSDE in the space of bounded processes (see for
%instance the recent papers by Hu, Imkeller and Muller
%and have a stability result prevailing in a wider context. This will
%allow us to obtain some results about the existence of a solution
%for a quadratic BSDE outside of the standard framework, moving away
%from the bounded case to the case where the terminal condition has
%%exponential moment.
%=========================================================
%s5 #&#
\section{Existence result for quadratic BSDEs}\label{sec5}
%=========================================================
The question of existence of bounded solutions for the classical
quadratic BSDEs in Brownian framework has been solved by Kobylanski
\cite{Kobylanski00}, using an exponential transformation as to
come back to the standard framework of a coefficient with linear
growth. A detailed review of the literature including the\vadjust{\goodbreak}
comparison theorem and different applications may be found in El
Karoui, Hamad\`ene and Matoussi~\cite{Elkaroui-Hamadene-Matoussi}.
Most of the recent papers focusing on financial applications of
quadratic BSDEs consider the situation where the martingale $M_{\cdot
}$ is
BMO
%as this gives a well-known framework for the existence of a
%solution for the BSDE in the space of bounded processes
(see, e.g., the recent papers by Hu, Imkeller and Muller \cite
{Hu-Imkeller-Muller}, Ankirchner, Imkeller and Reis
\cite{Ankirchner-Imkeller-Reis},
%Ankirchner, Imkeller and Popier
\cite{Morlais2}, or the Ph.D. thesis of dos Reis~\cite{DosReis}). From
Theorem~\ref{thevarqestimates}, such a framework is equivalent to look
at bounded solutions.
Briand and Hu~\cite{Briand-Hu} have been the first to extend
the previous results to unbounded solutions.
In all these papers, as in Kobylanski~\cite{Kobylanski00}, the main
difficulty is however to prove the strong convergence of the martingale
part.

The stability result we have
obtained in the previous section opens a new possible direction to
tackle this question. The idea is to approximate monotonically the
coefficient itself by coefficients with a linear and quadratic
growth, for which there are some results on the existence of
solution but also for which it is possible to take the limit thanks
to the stability Theorem~\ref{thmonotonecv}.
In our approach, we do not need this BMO framework
and have a stability result prevailing in a wider context,
%Existence results outside of the standard framework can be obtained,
moving away from the bounded case to the case where the terminal
condition has exponential moment. Indeed, having bounded
solutions is naturally replaced by belonging to the class ${\mathcal
S_Q}(|\eta_T|,\Lambda, C)$ as in the previous section, which reduces
to an exponential moment condition for $|\eta_T|$, when
$\Lambda\mbox{ and } C\equiv0.$ Recall that this last condition is
equivalent to have the absolute value of the solution in
the class ($\mathcal D_{\exp}$) when the coefficient does not depend on
$y$ [and $g(t,0,0) \equiv0$].

We start this section by looking more closely at the interrelationship between
%general
quadratic BSDEs and quadratic semimartingales,
%especially in the case
when the quadratic structure condition is saturated.
%quadratic BSDE}.

%that we recall here for the sake of clarity of the
%exposure: the various involved processes are some continuous
%predictable increasing process $K_{\cdot}$, two predictable bounded
%processes $l_{\cdot}$ and $c_{\cdot}$ such that $d\Lambda_t=l_tdK_t$
%and
%$dC_t=c_t dK_t$ and a $d$-dimensional continuous orthogonal
%martingale $N_{\cdot}$, for which the diagonal predictable quadratic
%variation matrix is dominated by $K_{\cdot}$, with
%$d\varq{N^i_t}=\gamma^{i}_t dK_t$ and $d\varq{N^i_t,N^j_t}=0$, if $
%i\neq j$. Finally $g(t,y,z)$ is a predictable process depending on
%$(y,z)\in\mathbb R \times\mathbb R^d$ in a continuous way.
%{\cgreen Pauline: je ne suis pas sure qu'il faille mettre cela ici, ou
%peut-etre comme remarque}\\
%Note that in this general framework, when $g$ has a linear growth
%(i.e. $|g(t,y,z)|
%z|$), existence results for Lipschitz coefficients may be found in
%El Karoui and Huang~\cite{NEKHuang} and can be easily extended to
%the case of a continuous coefficient with linear growth, following
%the arguments of Lepeltier and San Martin~\cite{LepeltierSanMartin98}.
%%%%%%%%%%%%%%%%%%%%%%%%%%%%%%%%%%%%%%
%s5.1 #&#
\subsection{\texorpdfstring{A canonical example: $q_\delta$-BSDE and entropic process}
{A canonical example: q delta-BSDE and entropic process}}\label{parq-BSDEs}
%%%%%%%%%%%%%%%%%%%%%%%%%%%%%%%%%%%%%%%
%As to illustrate how general quadratic BSDEs and quadratic
%semimartingales
%are interrelated,
We are focusing on simplest quadratic BSDEs when the structure
condition is saturated and the coefficient is simply denoted by
$q_{\delta}$. This framework has a particular importance in finance as
it corresponds to that of indifference pricing in incomplete markets
when using an exponential utility criterion (in general, in the bounded
case) as in Rouge and El Karoui~\cite{Rouge-ElKaroui}, and many other
papers (see, e.g., Mania and Schweizer~\cite{Mania-Schweizer})
or the recent book on indifference pricing edited by Carmona \cite
{Carmona}.

In this simple framework, it is interesting to consider the various
possible points of view. In particular, note that the two following
problems coincide in a Brownian framework:
\begin{longlist}[(ii)]
\item[(i)] First, finding a quadratic $q_\delta$-semimartingale
$Y_t=Y_0+M_t-\frac{\delta}{2}\varqM_t$ with terminal condition
$Y_T=\xi_T$.
%or a $\lu q_\delta$-semimartingale $Y_t=Y_0+M_t+\frac{\delta}{2}
We refer to the solution as a GBSDE($q_\delta, \xi_T)$-solution, where
$G$ stands for ``generalized.'' The process $-Y$ is a GBSDE-solution
associated with ($\lu q_\delta, -\xi_T)$.

\item[(ii)] In the second case, corresponding to the BSDE general framework
(Definition~\ref{defgeneralquadraticBSDE}), the problem is to find\vadjust{\goodbreak}
$(Y_{\cdot}, M_{\cdot} \equiv Z\ast N_{\cdot}+M^\bot_{\cdot})$,
such that $dY_t = -\frac
{\delta}{2}|\sqrt{\gamma_t}Z_t|^2 \,dK_t -Z_t \,dN_t-dM^\bot_t$ with
terminal condition $Y_T=\xi_T$. The similar equation with the opposite
process will be also considered. In the following, we refer to this
situation as $q$-BSDE.
\end{longlist}

%Note that $\delta$ is simply a scaling factor. Therefore, we make the
%normalizations $\delta= 1$ or $\delta=-1$ without any loss of
%generality, and denote $q \equiv q_1$, $\lu q \equiv q_{-1}$. Observe
%that the function $q$ is a convex function, while the function $\lu q$
%is concave.\\
%{\cgreen Pauline: doit-on preciser ici les notation $K$, $M^\bot$ et
%$N$?}.\\
%or for the $\lu q_\delta$-case,
%$dY_t = \frac{\delta}{2}|Z_t|^2 \,dK_t -Z_t.dN_t-dM^\bot_t $ with the
%same terminal condition $Y_T=\xi_T$.\\
Based on the previous results, we will consider these two questions in
parallel in the paragraphs below.
%The $q_\delta$ \rmi problem has been intensively study in Section 2
%and 3
%
%In this case, the coefficients are the functions
%$g(t,z)=q_{\delta}(z) \equiv\frac{\delta}{2} |z|^2$, and the
%associated BSDEs $-dY_t = \frac{\delta}{2}|Z_t|^2 \,dt -Z_t.dW_t,
% Y_T=\xi_T$ are the backward formulation of the
%$q_{\delta}$-semimartingales introduced in Section
%scaling factor, we make the normalizations $\delta= 1$ or $\delta
%=-1$, and denote $q=q_1, \lu q=q_{-1}$. In particular, we will see
%the essential role played by the variational approach, which will be
%generalized in the next subsection when applied to general quadratic
%BSDEs.
%%%%%%%%%%%%%%
%pa5.1.0.1 #&#
\subsubsection*{Summary of previous results on GBSDEs}
%%%%%%%%%%%%%%%%
The entropic process $\rho_{t}(\xi_T)$ defined earlier in
equation (\ref{eqentropicprocess}) as
$\ln\mathbb{E}[\exp(\xi_T)|\cF_t] \equiv\rho_{t}(\xi_T)$
appears naturally when studying such $(q, \mbox{ or } \lu q)$-GBSDEs.
Indeed, as presented in
the following proposition, if the terminal condition $\xi_T\in
\mathbb{L}_{\exp}^1$, then $\rho_{\cdot}(\xi_T)$ is a
$(\cD_{\exp})$-solution of $q$-GBSDE. The stronger assumption on the
terminal condition $|\xi_T|\in\mathbb{L}_{\exp}^1$ is used for the
estimates of the quadratic variation or for some stability result.
%
%pr5.1 #&#
%
\begin{eproposition}\label{propminimalentropic}
\begin{enumerate}[(iii)]
\item[(i)] Assume that $\xi_T \in\mathbb{L}_{\exp}^1$. Then
the entropic process $\rho_{\cdot}(\xi_T)$ is the unique $(\cD
_{\exp
})$-solution of the quadratic $\operatorname{GBSDE}(q,\xi_T)$, that is, there
exists a martingale $M^\rho_{\cdot}\in{\mathcal U}_{\exp}$ such
that
\[
d\rho_t(\xi_T)=-\tfrac{1}{2}\,d\varq{M^\rho}_t +dM^\rho_t,\qquad \rho
_T(\xi
_T)=\xi_T.
\]

Moreover, $\rho_{\cdot}(\xi_T)$ is minimal in the class of solutions
$Y_{\cdot}$:
$\rho_{\cdot}(\xi_T)\leq Y_{\cdot}$.
\item[(ii)]Assume that $-\xi_T \in\mathbb{L}_{\exp}^1$. The negative
entropic process $\urho_{\cdot}(\xi_T)$ is
{a} solution of the
%quadratic
$\operatorname{GBSDE}(\uq,\xi_T)$, i.e., there exists a martingale
$\underline{M}^\rho$ such that
\[
d\urho_t(\xi_T)=\tfrac{1}{2}\,d\varq{\underline{M}^{\rho}}_t
+d\underline
{M}^{\rho}_t,\qquad \urho_T(\xi_T)=\xi_T,
\]
but in general $\urho_{\cdot}(\xi_T)$ is not a $(\cD_{\exp})$-solution.

\item[(iii)]When $|\xi_T| \in\mathbb{L}_{\exp}^1$, then:
\begin{enumerate}[(b)]
\item[(a)] $\urho_t(\xi_T)$ is the maximal solution of the $\operatorname
{GBSDE}(\uq
,\xi_T)$.
\item[(b)] The martingales $M^\rho$ and $\underline{M}^{\rho}$ are in
$\mathbb
{H}^2$ and if $\xi_T$ is
bounded, they are BMO-martingales.

\item[(c)] If in addition $|\xi_T|+ \ln(|\xi_T|)\in\mathbb{L}_{\exp
}^1$, the
r.v. $\max|\rho_{0,T}(\xi_T)|$ and $\max|\urho_{0,T}(\xi_T)|$
belong to
$\mathbb{L}_{\exp}^1$. {Moreover, the following variational
representation holds true:}
%e24 #&#
%
\begin{equation}\label{eqdualdynamicentropy}
\rho_t(\xi_T)=\operatorname{ess}\sup_{M^{\mathbb Q}}\bigl\{\mathbb E_{\mathbb{Q}}
\bigl(\xi_T-\tfrac{1}{2}\varq{M^{\mathbb Q}}_{t,T}\bigr)|\cF_{t} | \mathbb
E_{\mathbb{Q}}(\varqM^{\mathbb Q}_{t,T})<+\infty\bigr\}.
\end{equation}
\end{enumerate}
\end{enumerate}
\end{eproposition}

\begin{pf}
\textup{(i)} From Section~\ref{subsectioncharactexpoinequality} and as
$\rho_ .(\xi_T)=\rho_ 0(\xi_T)+ r_{\cdot}(M)$, $\rho_{\cdot}(\xi
_T)$ is the unique
$(\cD_{\exp})$-solution for the GBSDE$(q,\xi_T)$, and the smallest in
the class of the $q$-semimartingale with the same terminal value.\vspace*{-6pt}
\begin{longlist}[(iii)]
%%%%
\item[(ii)] Since $-\xi_T \in\mathbb{L}_{\exp}^1$, the process
$\rho_{\cdot}(-\xi_T)$ is well defined\vspace*{1pt} in $(\expD)$ and $-\rho
_{\cdot}(-\xi_T)$
is solution of the $\lu q$-GBSDE, but not in general in the
class $(\expD)$.\vadjust{\goodbreak}

\item[(iii)] Assume both variables $\xi_T$ and $-\xi_T$ in
$\mathbb{L}^1_{\exp}$. Using the convexity of $\rho$, its follows
that $0=\rho_{\cdot}(0)\leq\frac{1}{2}(\rho_{\cdot}(\xi_T)-\urho
_{\cdot}(\xi_T))$. Then,
$\rho_{\cdot}(\xi_T)\in(\expD)$ implies $\urho_{\cdot}(\xi
_T)\in(\expD)$.

The comparison with the other solutions is a simple consequence of
the fact that $-Y$ is a solution of GBSDE$(q,-\xi_T)$, and therefore
bigger than $\rho_{\cdot}(-\xi_T)=-\urho_{\cdot}(\xi_T)$.
The rest of (iii) is a straightforward consequence of Theorem
\ref{thevarqestimates}.\quad\qed
\end{longlist}
\noqed\end{pf}

For lack of space, we will not further develop the variational point of
view, but this approach can be extended to $q$-BSDEs, using in
particular approximations based on the solutions of convex BSDEs with
linear growth (see, e.g., El Karoui, Hamad\`ene and Matoussi
\cite{Elkaroui-Hamadene-Matoussi}).
%
%pa5.1.0.2 #&#
\subsubsection*{\texorpdfstring{$(q \mbox{ or } \uq)$-BSDEs}{$(q, or \uq)$-BSDEs}} The
question of the existence of solutions of the $(q \mbox{ or }
\uq)$-BSDEs is more delicate to tackle and does not admit explicit
representation. These difficulties also appear in the Brownian
framework when the vector martingale~$N$ is defined from a limited
number of components of the generating Brownian motion. Different
methods can be used, the first one is based on linear growth
approximating solutions, whilst the second one uses the convexity of
the coefficient and represents solutions as value function of some
optimization problems. We now develop the first point of view.

In this case, the approximation is based on the coefficients $q_n(z)
\equiv\frac{1}{2}(|z|^{2}-(z-n)^{+2})=\frac{1}{2}(|z|^{2}\ind_{\{
|z|\leq n\}}+(n|z|-\frac{1}{2}n^2)\ind_{\{|z|> n\}})$ with linear and
quadratic growth, increasing to $q(z)$ when $n$ goes to infinity. For
$\xi_T\in{\mathbb L}^2$, using by the classical theory, the $\operatorname{BSDE}(q_n,
\xi_T)$ has a unique solution in ${\mathcal S}^2$, bounded if $\xi_T$
is bounded.
%
%pr5.2 #&#
%
\begin{eproposition}\label{propdualentropy}
%Let $q_n(z)=\frac{1}{2}|z|^{2}\ind_{\{|z|\leq n\}}+(n|z|-\demi n^2)
Let $|\eta_T| \in\mathbb L^1_{\exp}$, and $(\xi^n_T)$ a sequence of
increasing r.v., bounded by $|\eta_T|$ and converging a.s. to $\xi_T$.
\begin{longlist}[(iii)]
\item[(i)]Denote by $(Y^n_{\cdot}, Z^n_{\cdot}, M_{\cdot}^{n,\bot}) \in
\mathbb H^2(\mathbb
R^+)\otimes\mathbb H^2(\mathbb R^n)$ the unique solution of the $\operatorname{BSDE}(q_n, \xi_T)$. The process $Y^n_{\cdot}$ is a $\mathcal Q$-semimartingale
satisfying the entropic inequality $|Y^n_{\cdot}|\leq\rho_{\cdot
}(|\eta_T|)$.

\item[(ii)]The sequence $(q^n, Y^n_{\cdot}, Z^n_{\cdot}, M_{\cdot
}^{n,\bot}) $ satisfies
the hypothesis of Theorem~\ref{thmonotonecv} and strongly converges to
$(Y_{\cdot}, Z_{\cdot}, M_{\cdot}^{\bot}) $, minimal solution of\break
$\operatorname{BSDE}(q,\xi_T)$ such
that $|Y_{\cdot}|\leq\rho_{\cdot}(|\eta_T|)$, with the variational
representation
%e25 #&#
%
\begin{equation}\label{eqdualdynamicq-entropy}
Y_t(\xi_T)=\operatorname{ess}\sup_{\nu}\biggl\{\mathbb E_{\mathbb{Q^\nu}}\biggl(\xi
_T-\frac{1}{2}\int_t^T\bigl|\sqrt{\gamma_s }\nu_s\bigr|^2\,dK_s\big|\cF_{t} \biggr)\biggr\},
\end{equation}
where $\mathbb{Q^\nu}$ is the probability with density $\mathcal
E_{\cdot}(\nu
\star N)$ with finite entropy\break
$\mathbb E_{\mathbb{Q^\nu}}( \int_0^T|\sqrt{\gamma_t }\nu
_t|^2\,dK_t)<+\infty$.

\item[(iii)] Uniqueness holds in the class of solutions $Y$ such that
$|Y_{\cdot}|\leq\rho_{\cdot}(|\xi_T|)$, and $|\xi_T|$ is bounded
or such that $\rho
_\delta(|\xi_T|)$ for any $\delta>0$.
\end{longlist}
\end{eproposition}

\begin{pf} (i) Its is clear that $Y^n_{\cdot}$ is a $\mathcal
Q$-semimartingale, bounded if $\xi_T$ is bounded. Then $|Y^n_{\cdot}| $
belongs to the class $({\mathcal D}_{\exp})$ and satisfies the entropic
inequality $|Y^n_{\cdot}|\leq\rho_{\cdot}(|\xi_T|)$. Since both processes
$Y^n_{\cdot}(\xi_T)$ and $\rho_{\cdot}(\xi_T)$ are monotone with
respect to their
terminal condition (by approximating~$\xi_T)$ by bounded random
variables, the entropic inequality holds at the limit under the
assumption $|\xi_T|\in\mathbb L^1_{\exp}$.\vspace*{-6pt}
\begin{longlist}[(iii)]
\item[(ii)] The first result is a direct consequence of the stability result
given in Theorem~\ref{thmonotonecv}.
The variational representation for $Y^n$ is as in \eqref
{eqdualdynamicq-entropy} with the restriction that $\nu$ is bounded by
$n$. That is a standard result on convex BSDEs with uniformly linear
growth (see, e.g., El Karoui, Peng and Quenez~\cite
{ElKaroui-Peng-Quenez} or Theorem 8.7
in El Karoui, Hamad\`{e}ne and Matoussi \cite
{Elkaroui-Hamadene-Matoussi}). Thanks to entropy result in Section
\ref{changeprobaentropy}, the representation \eqref
{eqdualdynamicq-entropy} pass to the limit, since $\xi_T$ is $\mathbb
Q^\nu$-integrable.

\item[(iii)] Let $Y$ be a solution satisfying $|Y_{\cdot}|\leq\rho_{\cdot
}(|\xi_T|)$. We
first assume that $|\xi_T|$ is bounded, so that all solutions are
bounded and the associated martingales $M_{\cdot}$, $M^{\bot}_{\cdot
}$ and
$M_{\cdot}-M^{\bot}_{\cdot}$ are BMO-martingales.

Denote by $Y^i_{\cdot}$ and $Y^i_{\cdot}$ two solutions satisfying
the entropic
inequalities with two bounded terminal conditions. Using the same
notation than in the proof of Theorem~\ref{thNicolas}, we observe that
the difference $Y^{i,j}_{\cdot} \equiv Y^i_{\cdot}-Y^j_{\cdot}$
verifies a linear BSDE,
with linear growth condition with respect to another probability measure,
\begin{eqnarray*}
dY^{i,j}_t&=&-\tfrac{1}{2}\bigl(\bigl|\sqrt{\gamma_t }Z^i_t\bigr|^2-\bigl|\sqrt{\gamma_t
}Z^j_t\bigr|^2\bigr)\,dK_t+dM^{i,j}_t\\
&=&-\tfrac{1}{2}\bigl(\sqrt{\gamma_t }Z^{i,j}_t.\sqrt{\gamma_t
}(Z^i_t+Z^j_t)\bigr)\,dK_t+Z^{i,j}_t.\,dN_t+dM^{i,j,\bot}_t\\
%&=-\demi(|\sqrt{\gamma_t }Z'_t|^2-|\sqrt{\gamma_t }Z_t|^2-2
&=&Z^{i,j}_t \bigl(dN_t-\tfrac{1}{2}\sqrt{\gamma_t
}(Z^i_t+Z^j_t)\,dK_t\bigr)+dM^{i,j,\bot}_t.
\end{eqnarray*}
Since $Y^i_{\cdot}$ and $Y^j_{\cdot}$ are bounded solutions, by
Theorem \ref
{thevarqestimates}, the martingales $M^i_{\cdot}$, and $M^j_{\cdot}$ are
BMO-martingales, implying that the quadratic variation of $\frac{1}{2}
(Z^i_{\cdot}+Z^j_{\cdot})\star N_{\cdot}$ is also conditionally
bounded, and then $\frac{1}{2}
(Z^1_{\cdot}+Z^2_{\cdot})\star N_{\cdot}$ is a BMO-martingale. By
Girsanov theorem,
$\mathcal E(\frac{1}{2}(Z^1_{\cdot}+Z^2_{\cdot})\star N)_{\cdot}$
is a
u.i. exponential
martingale defining a new probability measure $\mathbb Q^{(i+j)}$ such
that $dN^{(i+j)}_t\equiv dN_t-\frac{1}{2}\sqrt(Z^i_t+Z^j_t)\,dK_t)$
is a
$\mathbb Q^{(i+j)}$-local martingale with the same quadratic variation
as $N_{\cdot}$. Moreover, $M^{i,j,\bot}_{\cdot}$ is still a $\mathbb
Q^{(i+j)}$-local martingale, orthogonal to $N^{(i+j)}_{\cdot}.$
Then, $Y^{i,j}_{\cdot} $ is a bounded $\mathbb Q^{(i+j)}$ local martingale,
and so a true martingale and $Y^{i,j}_{\cdot} =\mathbb E_{\mathbb
{Q}^{(i+j)}}(Y^{i,j}_T|\cF_{\cdot})$. Uniqueness and comparison
theorem are
easily deduced of this property.

In the general case, the difficulty is to show directly that $\mathcal
E(\frac{1}{2}(Z^1_{\cdot}+\break Z^2_{\cdot})\star N)_{\cdot}$ is u.i.
martingale, given that
$\mathcal E(M^i)_{\cdot}$ and $\mathcal E(M^j)_{\cdot}$ are uniformly
integrable.

Under the assumptions of exponential moments of any order, uniqueness
has been proved first by Briand and Hu~\cite{Briand-Huconvex} and Mocha
and Westray~\cite{Mocha-Westray}.\quad\qed
\end{longlist}
\noqed\end{pf}

\subsection{\texorpdfstring{Existence result for BSDEs in the class ${\mathcal S_Q}(|\eta_T|,\Lambda, C)$}
{Existence result for BSDEs in the class S Q(|eta T|, Lambda, C)}}
%=================================================
We are now interested in quadratic BSDEs satisfying the general
structure condition $| g(\cdot,t,y,z)|\leq\kappa(t,y,z) \equiv
|l_{t}|+c_t|y|+\frac{1}{2}|z|^{2}, d\mathbb{P}\otimes dK_{\cdot
}\ \mathrm{a.s.}$, and are looking for solution in the class ${\mathcal
S_Q}(|\eta
_T|,\Lambda, C)$ only. As before, the method
relies on a regularization of the
quadratic coefficient it-self through inf-convolution as to
transform it into a coefficient with \textit{both} linear and quadratic
growth. This double structure of the transformed
coefficient leads to results both in terms of existence and
estimation.
%Indeed, in the light of the previous results of this
%paper, we can consider these BSDEs as quadratic semimartingales, and
%so it seems very natural to study them under the previous Hypothesis
%condition.
The previous stability Theorem~\ref{thNicolas} can then
be applied to obtain the existence of a solution, after having
proved that
the approximate solutions are also ${\mathcal S_Q}(|\eta_T|,\Lambda,
C)$-semimartingales.
%
%pa5.2.0.1 #&#
\subsubsection*{Regularization of the coefficient through inf-convolution}
%
%The functions $q_n(z)$
%introduced in Proposition~\ref{propdualentropy} are an example of
%regularization by inf-convolution of the canonical function $q(z)=
The proof of this fundamental result is based on the following lemma
involving classical regularization by inf-convolution techniques
introduced by Lepeltier and San Martin~\cite{LepeltierSanMartin97}
in a BSDEs framework. Let us first observe that the appropriate
regularization when dealing with $\uq(z)=-\frac{1}{2}|z|^2$ is a
sup-convolution since $\uq(z)$ is concave. To overcome this difficulty,
we proceed in two steps, by first assuming that $g$ is bounded from
below by some basic function with both a linear and quadratic growth
$\lu{\kappa}_p$, where $-\lu{\kappa}_p(t,y,z) =\kappa
_p(t,y,z)\equiv
l_t+c_t |y|+ q_p(z)$ with $q_p(z)=\frac{1}{2}|z|^{2}\ind_{\{|z|\leq
p\}
}+(p|z|-\frac{1}{2}p^2)\ind_{\{|z|> p\}}$. When $p=1$, $\kappa_1(t,y,z)
\equiv\kappa(t,y,z)=l_t+c_t |y|+ q(z)$ with $q(z)=\frac{1}{2}|z|^2$.
%
%le5.3 #&#
%
\begin{elemma}\label{leminfconvol}
Let $g:{\mathbb R}\times{\mathbb R}^n\rightarrow{\mathbb
R}$ be a continuous function with linear growth in y, and quadratic
growth in $z$, bounded from below by some function $\lu{\kappa
}_p(t,y,z)=-(l_t+c_t |y|+ q_p(z))$ and from above by $\kappa(t,y,z)$:
%e26 #&#
%
\begin{equation}
\label{eqgconstraint}
\lu{\kappa}_p(t,y,z) \leq g(t,y,z) \leq\kappa(t,y,z),
\end{equation}
where the processes $c_{\cdot}$ and $l_{\cdot}$ are bounded by some universal
constant $\bar{C}$.

The regularizing functions are the convex functions with
linear growth $b_n(u,w)\equiv n |u|+n |w|$.The sequences $\lu{\kappa
}_{n,p}$, $\kappa_n$ and $g_n$ are defined, respectively, as the
inf-convolution
of the functions $\lu{\kappa}_p$, $\kappa$ and $g$ with the function~$b_n$,
\begin{eqnarray*}
\lu{\kappa}_{n,p}(t,y,z)&=&\lu{\kappa}_p \,\square\, b_n(t,y,z),\qquad
\kappa_n(t,y,z)=\kappa\,\square\, b_n(t,y,z),\\
g_n(t,y,z)&=&g \,\square\, b_n(t,y,z)=\inf_{u,w }
\bigl(g(t,u,w)+n|y-u|+n|z-w|\bigr)
\end{eqnarray*}
have the following properties, for $n\geq\sup(\bar{C},p)$:
\begin{longlist}[(iii)]
\item[(i)] $ \kappa_n(t,y,z)\,{=}\,l_t+c_t|y|+q_n(z)\leq l_t+c_t|y|+ \frac{1}{2}
|z|^2, \lu{\kappa}_{n,p}(t,y,z)\,{=}\, \lu{\kappa}_{p}(t,y,z)$;

\item[(ii)] $|g_n(t,y,z)|\leq l_t+c_t|y|+\sup(q_p(z), q_n(z))=\kappa
_n(t,y,z)\leq l_t+c_t |y|+ \frac{1}{2}|z|^{2} $;

\item[(iii)]the sequences $g_n$ and $\kappa_n$ are increasing;

\item[(iv)] the Lipschitz constant of functions $g_n$ is $n$;

\item[(v)] if $(y_n,z_n)\rightarrow(y,z)$, then $g_n(t,y_n, z_n)\rightarrow
g(t,y,z)$.\vadjust{\goodbreak}
\end{longlist}
\end{elemma}

In this lemma, the various functions are regularized through the
Lipschitzian regularization, whist the function $\kappa_n$ is the
Moreau--Yoshida regularization of $b_n$ (see Hiriart-Urruty and Lemar\'
{e}chal~\cite{Hiriart-Urruty}, Chapter E, for more details).
%The inf-convolution of coefficients have been extensively studied in
%the context of risk transfer by Barrieu and El Karoui
%sharing using dynamic financial risk measures as criterion can be
%reduced to the study of the inf-convolution of the associated
%coefficients.
%
%pa5.2.0.2 #&#
\subsubsection*{Existence result}
The important point now is to prove that the solutions to the
$\operatorname{BSDEs}(g_n, \xi_T)$ which are Lipschitz with linear growth, are
in the class ${\mathcal S_Q}(|\eta_T|,\Lambda, C)$ when $\E[\exp
(\bar
X^{\Lambda,C}_T(|\eta_T|))]<+\infty$.
%The argument is a slight extension of the one used in
%the proof of
%Proposition~\ref{propdualentropy}.
%le5.4 #&#
%
\begin{elemma}\label{lemlinearquadraticgrowth}
%Let $|\eta_T|$ be a $\cF_T$-r.v. such that $\E[\exp(\bar X^{
Let $g$ and $g_n$, $\kappa$ and $\kappa_n$ as in Lemma \ref
{leminfconvol}. The coefficients $g_n$ and $\kappa_n$ are standard
uniformly Lipschitz coefficients.
For any $|\xi_T|\leq|\eta_T|$, let $(Y^n_{\cdot},Z^n_{\cdot
},M^{n,\bot}_{\cdot})$ and
$(U^n_{\cdot},V^n_{\cdot},W^{n,\bot}_{\cdot})$ be the unique
solution of the $\operatorname{BSDE}(g_n,\xi
_T)$ and $\operatorname{BSDE}(\kappa_n,|\eta_T|)$ in the appropriate space.
\begin{longlist}[(ii)]
\item[(i)] The sequences $(Y^n_{\cdot})$ and $(U^n_{\cdot})$ are
increasing, and satisfy
the entropic inequality, $|Y^n_{\cdot}|\leq U^n_{\cdot} \leq\rho
_{\cdot}(e^{C_{\cdot,T} }
|\eta_T|+\int_{\cdot}^T
e^{C_{\cdot,t}}\,d\Lambda_t ), a.s.$

%$ and $U^n_{\cdot}$ satisfies the entropic inequality $U^n_{\cdot}\leq
%e^{C_{\cdot,t}}\,d\Lambda_t ), a.s.$\\
%U^n_\sigma\equiv\operatorname{ess}\sup_{ |\nu_t|\leq n}
%e^{C_{\sigma,t}}\,d\Lambda_t - \demi\int_\sigma^T
%e^{C_{\sigma,t}}\,d\varq{\nu.N}_t|\cF_\sigma] ;
%The random variable $U^n_\sigma$ is dominated by
%$\rho_\sigma(e^{C_{\sigma,T} } |\eta_T|+\int_\sigma^T
%e^{C_{\sigma,t}}\,d\Lambda_t)$, and
Both sequences $(Y^n_{\cdot})$ and $(U^n_{\cdot})$
are ${\mathcal S_Q}(|\eta_T|,\Lambda, C)$-quadratic semimartingales.

\item[(ii)] The sequence $(Y^n_{\cdot},Z^n_{\cdot},M^{n,\bot}_{\cdot})$
converges uniformly in
probability to a minimal solution $(Y_{\cdot},Z_{\cdot}, M^{\bot
}_{\cdot})$ of the
$\operatorname{BSDE}(g, \xi_T)$.
\end{longlist}
\end{elemma}

\begin{pf} The proof relies on classical properties of BSDEs solutions
associated with standard coefficients (with linear growth), in a
$\mathbb H^2$-space. In particular, existence, uniqueness and
comparison hold true in this case, that implies (i).
\begin{longlist}[(ii)]
\item[(i)] First, assume that ${\bar X}^{\Lambda,C}_T$ is bounded. The
solutions $U^n$ are bounded and the entropic inequality is valid.
%The submartingale $\exp({\bar X}^{\Lambda,C}_{\cdot}(U^n))$ is also a
%bounded, and then the process $U^n$ satisfies the desired inequalities.
Since these inequalities are stable when taking increasing limit with
respect to $\Lambda,C, \eta$, the same inequalities hold still true
under the assumption ${\bar X}^{\Lambda,C}_T(\eta_T)\in\mathbb
L^1_{\exp}$.
%
%place the bounded processes $C$ and $\Lambda$ by the processes $
Then, by construction, $(Y^n_{\cdot})$ and $(U^n_{\cdot})$ are
${\mathcal S_Q}(|\eta_T|,\Lambda, C)$- quadratic semimartingales.

\item[(ii)] Finally, using Theorem~\ref{thNicolas}, we
obtain the convergence of this sequence to a solution of the
$\operatorname{BSDE}(g, \xi_T)$ in the space ${\mathcal
S_Q}(|\eta_T|,\Lambda,C)$.\quad\qed
\end{longlist}
\noqed\end{pf}

It remains to overcome the assumption made on the coefficient of a
linear quadratic growth lower bound.
Given a coefficient $g$ with decomposition $g=g^+-g^-$, where both
positive functions $g^+$ and $g^-$ have the same quadratic structure.
Let $g_p \equiv g^+-g^-\,\square\, b_p$. Then $ g_p$ satisfies Condition
(\ref{eqgconstraint}), and the $\operatorname{BSDE}(g_p, \xi_T)$ admits a minimal
solution; the sequence of solutions $Y^p$ is decreasing, and belongs to
the space ${\mathcal S_Q}(|\eta_T|,\Lambda,C)$. Once again, we use the
stability theorem to conclude that the sequence $Y^p$ converges to a
solution of the $\operatorname{BSDE}(g,\xi_T)$. We summarize the general form of our
results in the following theorem.

%
%th5.5 #&#
%
\begin{theorem}\label{lemgeneralexistenceresult} Let us consider a
general $\operatorname{BSDE}(g,\xi_T)$,
where $\xi_T$ be a $\cF_T$-random variable such that
$\E[\exp(e^{C_{T}} | \delta\xi_T|+ \int_0^Te^{C_s}\,d\Lambda
_s)]<+\infty$.

The coefficient $g(t,y,z)$ is satisfying
the quadratic structure condition (\ref{eqquadraticgrowth}), $|
g(\cdot,t,y,z)|\leq\frac{1}{\delta} |l_{t}|+c_t|y|+\frac{\delta
}{2}|z|^{2}$.

Then, there exists at least a solution $(Y,Z,M^{\bot})$ in ${\mathcal
S_Q}(|\xi_T|,\Lambda, C,\delta)$ of the $\operatorname{BSDE}(g,\xi_T)$.
\end{theorem}

%re5 #&#
%
\begin{remark} When both $\Lambda, C \equiv0$, as in the
framework of cash additive risk measures, the theorem simply states: if
$| g(\cdot,t,y,z)|\leq\frac{\delta}{2}|z|^{2}$, and $\mathbb E[\exp
(\delta|\xi_T|)]<+\infty$, their exists at least a solution in the
class $({\mathcal D}_{\exp}$).
\end{remark}

%and~\cite{Delbaen-Hu-Richou}, %Bao, Delbaen and Hu and Delbaen, Hu and
%%Richou have shown} that when the coefficient $g$ is %convex, we can
%simply assume that $\xi_T\in\mathbb{L}^1_{exp}$ and $g(t,y,z)\leq l+

%A straightforward corollary of this theorem is given below as:
%$|g(t,y,z)|\leq|%g_0(t)|+c_t |y|+\gamma_t |z|$. Let $l_t=|g_0(t)|+
%There exists at least a solution of the linear growth BSDE in the
%space $%{\mathcal
%S_Q}(|\xi_T|,\Lambda, C)$ for any $\xi_T$ such that $\E[
%%%%%%%%%%%%%%%%%%%
%pa5.2.0.3 #&#
\subsubsection*{Comment on the uniqueness of the solution}
%%%%%%%%%%%%%%%%%%%%%%%%%%%%%%%%
The question of the uniqueness of the solution to a general
quadratic BSDE is not trivial. In the standard framework where
the terminal condition is bounded, Kobylanski~\cite{Kobylanski00}
obtains the uniqueness of the solution under some Lipschitz style
assumptions. Recently, Tevzadze~\cite{Tevzadze} gives a direct proof
of uniqueness still in the bounded case. In the case of an unbounded
terminal condition, Briand and Hu~\cite{Briand-Huconvex} work under
the additional assumption that the coefficient $g$ is convex with
respect to the variable $z$. This allows them to derive the
comparison theorem, which is needed to obtain the uniqueness. Their
methodology can be adapted and generalized to our framework without
any particular difficulty. In a very recent paper
\cite{Mocha-Westray}, Mocha and Westray have considered general
quadratic BSDEs under some stronger assumptions of exponential
moment of order $p>1$ and boundedness of the increasing processes.
They obtain some interesting results for the uniqueness of the
solution. The convex case has been also studied in Delbaen, Hu and Richou
\cite{Delbaen-Hu-Richou}
under weaker assumptions.

In this paper, we study the stability and convergence of some general
quadratic semimartingales. The general stability result (Theorem \ref
{thNicolas}), including the strong convergence of the martingale parts
in various spaces ranging from $\mathbb{H}^1$ to BMO, is derived under
some mild integrability condition on the exponential of the terminal
value of the semimartingale. This strong convergence result is then
used to prove the existence of solutions of general quadratic BSDEs
under minimal exponential integrability assumptions, relying on a
regularization in both linear-quadratic growth of the quadratic
coefficient itself. On the contrary to most of the existing literature,
it does not involve the seminal result of Kobylanski \cite
{Kobylanski00} on bounded solutions. As previously mentioned, this
approach has also other potential applications such as numerical
simulations of quadratic BSDEs, study in terms of risk measures and
dual representation, solving of associated HJB-type equations\ldots. The
various results obtained in the paper can be extended to jump processes.

%==========================================================
\section*{Acknowledgments}
%===========================================================
Both authors would like to thank Mingyu Xu and Anis Matoussi for their
helpful comments and discussions at various stages of the paper. The
paper has also greatly benefited from a careful reading and helpful
suggestions from two anonymous referees. Finally, the authors would
like to give some special thanks to Arthur who waited patiently before arriving.

%=============================================================
% imsref loaded by akundreckaite, 2012-06-04 08:39:51
% imsref loaded by akundreckaite, 2012-06-04 13:20:50
% imsref loaded by akundreckaite, 2012-06-04 13:23:42
% imsref loaded by akundreckaite, 2012-06-04 13:35:12
%

%

%suskaldyti doi

\printaddresses


\begin{thebibliography}{46}
% BibTex style file: ims.bst, 2011-05-30
% Default style options (sort=0,type=number).
% Used options (sort=1,type=number).

%b1 #&#
\bibitem{Ankirchner-Imkeller-Reis}
%
\begin{barticle}[mr]
\bauthor{\bsnm{Ankirchner},~\bfnm{Stefan}\binits{S.}},
\bauthor{\bsnm{Imkeller},~\bfnm{Peter}\binits{P.}} \AND
\bauthor{\bsnm{Dos~Reis},~\bfnm{Gon{\c{c}}alo}\binits{G.}}
(\byear{2007}).
\btitle{Classical and variational differentiability of {BSDE}s with quadratic
growth}.
\bjournal{Electron. J. Probab.}
\bvolume{12}
\bpages{1418--1453 (electronic)}.
\bid{doi={10.1214/EJP.v12-462}, issn={1083-6489}, mr={2354164}}
\bptok{imsref}%
\end{barticle}
%
\endbibitem

%b2 #&#
\bibitem{Azema-Gundy-Yor}
%
\begin{bincollection}[mr]
\bauthor{\bsnm{Az{\'e}ma},~\bfnm{J.}\binits{J.}},
\bauthor{\bsnm{Gundy},~\bfnm{R.~F.}\binits{R.~F.}} \AND
\bauthor{\bsnm{Yor},~\bfnm{M.}\binits{M.}}
(\byear{1980}).
\btitle{Sur l'int\'egrabilit\'e uniforme des martingales continues}.
In \bbooktitle{Seminar on {P}robability, {XIV} (Paris, 1978/1979) ({F}rench)}.
\bseries{Lecture Notes in Math.}
\bvolume{784}
\bpages{53--61}.
\bpublisher{Springer}, \baddress{Berlin}.
\bid{mr={0580108}}
\bptok{imsref}%
\end{bincollection}
%
\endbibitem

%%b3 #&#
%%
%(\byear{2010}).
%In \bbooktitle{Contemporary Quantitative Finance}
%(\beditor{\binits{C.} \bsnm{Chiarella}}
%%

%b4 #&#
\bibitem{Barlow-Protter}
%
\begin{bincollection}[mr]
\bauthor{\bsnm{Barlow},~\bfnm{Martin~T.}\binits{M.~T.}} \AND
\bauthor{\bsnm{Protter},~\bfnm{Philip}\binits{P.}}
(\byear{1990}).
\btitle{On convergence of semimartingales}.
In \bbooktitle{S\'eminaire de {P}robabilit\'es, {XXIV}, 1988/89}.
\bseries{Lecture Notes in Math.}
\bvolume{1426}
\bpages{188--193}.
\bpublisher{Springer}, \baddress{Berlin}.
\bid{doi={10.1007/BFb0083765}, mr={1071540}}
\bptok{imsref}%
\end{bincollection}
%
\endbibitem

%b5 #&#
\bibitem{Barrieu-ElKaroui5}
%
\begin{bincollection}[mr]
\bauthor{\bsnm{Barrieu},~\bfnm{Pauline}\binits{P.}} \AND
\bauthor{\bsnm{El~Karoui},~\bfnm{Nicole}\binits{N.}}
(\byear{2004}).
\btitle{Optimal derivatives design under dynamic risk measures}.
In \bbooktitle{Mathematics of Finance}
(\beditor{\binits{G.}~\bsnm{Yin}} \AND
\beditor{\binits{Q.}~\bsnm{Zhang}}, eds.).
\bseries{Contemp. Math.}
\bvolume{351}
\bpages{13--25}.
\bpublisher{Amer. Math. Soc.}, \baddress{Providence, RI}.
\bid{doi={10.1090/conm/351/06389}, mr={2076287}}
\bptok{imsref}%
\end{bincollection}
%
\endbibitem

%b6 #&#
\bibitem{Barrieu-ElKaroui6}
%
\begin{bincollection}[auto:STB|2012/06/01|15:11:24]
\bauthor{\bsnm{Barrieu},~\bfnm{P.}\binits{P.}} \AND
\bauthor{\bsnm{El~Karoui},~\bfnm{N.}\binits{N.}}
(\byear{2009}).
\btitle{Pricing, hedging and optimally designing derivatives via minimization
of risk measures}.
In \bbooktitle{Volume on Indifference Pricing}
(\beditor{\binits{R.}~\bsnm{Carmona}}, ed.)
\bpages{77--146}.
\bpublisher{Princeton Univ. Press}, \baddress{Princeton}.
\bptok{imsref}%
\end{bincollection}
%
\endbibitem

%b7 #&#
\bibitem{Briand-Hu}
%
\begin{barticle}[mr]
\bauthor{\bsnm{Briand},~\bfnm{Philippe}\binits{P.}} \AND
\bauthor{\bsnm{Hu},~\bfnm{Ying}\binits{Y.}}
(\byear{2006}).
\btitle{B{SDE} with quadratic growth and unbounded terminal value}.
\bjournal{Probab. Theory Related Fields}
\bvolume{136}
\bpages{604--618}.
\bid{doi={10.1007/s00440-006-0497-0}, issn={0178-8051}, mr={2257138}}
\bptok{imsref}%
\end{barticle}
%
\endbibitem

%b8 #&#
\bibitem{Briand-Huconvex}
%
\begin{barticle}[mr]
\bauthor{\bsnm{Briand},~\bfnm{Philippe}\binits{P.}} \AND
\bauthor{\bsnm{Hu},~\bfnm{Ying}\binits{Y.}}
(\byear{2008}).
\btitle{Quadratic {BSDE}s with convex generators and unbounded terminal
conditions}.
\bjournal{Probab. Theory Related Fields}
\bvolume{141}
\bpages{543--567}.
\bid{doi={10.1007/s00440-007-0093-y}, issn={0178-8051}, mr={2391164}}
\bptok{imsref}%
\end{barticle}
%
\endbibitem

%b9 #&#
\bibitem{Carmona}
%
\begin{bbook}[mr]
\beditor{\bsnm{Carmona},~\bfnm{Ren{\'e}}\binits{R.}}, ed.
(\byear{2009}).
\btitle{Indifference Pricing: Theory and Applications}.
\bpublisher{Princeton Univ. Press}, \baddress{Princeton, NJ}.
\bid{mr={2547456}}
\bptok{imsref}%
\end{bbook}
%
\endbibitem

%b10 #&#
\bibitem{ChoulliSchweizer}
%
\begin{bmisc}[auto:STB|2012/06/01|15:11:24]
\bauthor{\bsnm{Choulli},~\bfnm{T.}\binits{T.}} \AND
\bauthor{\bsnm{Schweizer},~\bfnm{M.}\binits{M.}}
(\byear{2011}).
\bhowpublished{Stability of Sigma-martingale densities in $L\log L$
under an
equivalent change of measure. NCCR FINRISK Working Paper 676, ETH Zurich}.
\bptok{imsref}%
\end{bmisc}
%
\endbibitem

%%b11 #&#
%%
%(\byear{2002}).
%%

%b12 #&#
\bibitem{Delbaen-Hu-Bao}
%
\begin{barticle}[mr]
\bauthor{\bsnm{Delbaen},~\bfnm{Freddy}\binits{F.}},
\bauthor{\bsnm{Hu},~\bfnm{Ying}\binits{Y.}} \AND
\bauthor{\bsnm{Bao},~\bfnm{Xiaobo}\binits{X.}}
(\byear{2011}).
\btitle{Backward {SDE}s with superquadratic growth}.
\bjournal{Probab. Theory Related Fields}
\bvolume{150}
\bpages{145--192}.
\bid{doi={10.1007/s00440-010-0271-1}, issn={0178-8051}, mr={2800907}}
\bptok{imsref}%
\end{barticle}
%
\endbibitem

%b13 #&#
\bibitem{Delbaen-Hu-Richou}
%
\begin{barticle}[mr]
\bauthor{\bsnm{Delbaen},~\bfnm{Freddy}\binits{F.}},
\bauthor{\bsnm{Hu},~\bfnm{Ying}\binits{Y.}} \AND
\bauthor{\bsnm{Richou},~\bfnm{Adrien}\binits{A.}}
(\byear{2011}).
\btitle{On the uniqueness of solutions to quadratic {BSDE}s with convex
generators and unbounded terminal conditions}.
\bjournal{Ann. Inst. Henri Poincar\'e Probab. Stat.}
\bvolume{47}
\bpages{559--574}.
\bid{doi={10.1214/10-AIHP372}, issn={0246-0203}, mr={2814423}}
\bptok{imsref}%
\end{barticle}
%
\endbibitem

%b14 #&#
\bibitem{Dellacherie}
%
\begin{bincollection}[mr]
\bauthor{\bsnm{Dellacherie},~\bfnm{C.}\binits{C.}}
(\byear{1979}).
\btitle{In\'egalit\'es de convexit\'e pour les processus croissants et les
sousmartingales}.
In \bbooktitle{S\'eminaire de {P}robabilit\'es, {XIII} ({U}niv. {S}trasbourg,
{S}trasbourg, 1977/78)}.
\bseries{Lecture Notes in Math.}
\bvolume{721}
\bpages{371--377}.
\bpublisher{Springer}, \baddress{Berlin}.
\bid{mr={0544807}}
\bptok{imsref}%
\end{bincollection}
%
\endbibitem

%b15 #&#
\bibitem{Dellacherie-Meyer}
%
\begin{bbook}[mr]
\bauthor{\bsnm{Dellacherie},~\bfnm{Claude}\binits{C.}} \AND
\bauthor{\bsnm{Meyer},~\bfnm{Paul-Andr{\'e}}\binits{P.-A.}}
(\byear{1975}).
\btitle{Probabilit\'es et Potentiel}.
\bpublisher{Hermann}, \baddress{Paris}.
\bid{mr={0488194}}
\bptok{imsref}%
\end{bbook}
%
\endbibitem

%b16 #&#
\bibitem{DosReis}
%
\begin{bbook}[auto:STB|2012/06/01|15:11:24]
\bauthor{\bsnm{Dos~Reis},~\bfnm{G.}\binits{G.}}
(\byear{2011}).
\btitle{Some Advances on Quadratic BSDE: Theory--Numerics--Applications}.
\bpublisher{Lap Lambert Academic Publishing}, \baddress{Germany}.
\bptok{imsref}%
\end{bbook}
%
\endbibitem

%b17 #&#
\bibitem{ElKaroui-Hamadene}
%
\begin{barticle}[mr]
\bauthor{\bsnm{El-Karoui},~\bfnm{N.}\binits{N.}} \AND
\bauthor{\bsnm{Hamad{\`e}ne},~\bfnm{S.}\binits{S.}}
(\byear{2003}).
\btitle{B{SDE}s and risk-sensitive control, zero-sum and nonzero-sum game
problems of stochastic functional differential equations}.
\bjournal{Stochastic Process. Appl.}
\bvolume{107}
\bpages{145--169}.
\bid{doi={10.1016/S0304-4149(03)00059-0}, issn={0304-4149}, mr={1995925}}
\bptok{imsref}%
\end{barticle}
%
\endbibitem

%b18 #&#
\bibitem{Elkaroui-Hamadene-Matoussi}
%
\begin{bincollection}[auto:STB|2012/06/01|15:11:24]
\bauthor{\bsnm{El~Karoui},~\bfnm{N.}\binits{N.}},
\bauthor{\bsnm{Hamad{\`e}ne},~\bfnm{S.}\binits{S.}} \AND
\bauthor{\bsnm{Matoussi},~\bfnm{A.}\binits{A.}}
(\byear{2009}).
\btitle{BSDEs and applications}.
In \bbooktitle{Volume on Indifference Pricing}
(\beditor{\binits{R.}~\bsnm{Carmona}}, ed.)
\bpages{267--320}.
\bpublisher{Princeton Univ. Press}, \baddress{Princeton}.
\bptok{imsref}%
\end{bincollection}
%
\endbibitem

%b19 #&#
\bibitem{NEKHuang}
%
\begin{bincollection}[mr]
\bauthor{\bsnm{El~Karoui},~\bfnm{N.}\binits{N.}} \AND
\bauthor{\bsnm{Huang},~\bfnm{S.~J.}\binits{S.~J.}}
(\byear{1997}).
\btitle{A general result of existence and uniqueness of backward stochastic
differential equations}.
In \bbooktitle{Backward Stochastic Differential Equations ({P}aris,
1995--1996)}
(\beditor{\binits{N.}~\bsnm{El Karoui}} \AND
\beditor{\binits{L.}~\bsnm{Mazliak}},
eds.).
\bseries{Pitman Res. Notes Math. Ser.}
\bvolume{364}
\bpages{27--36}.
\bpublisher{Longman}, \baddress{Harlow}.
\bid{mr={1752673}}
\bptok{imsref}%
\end{bincollection}
%
\endbibitem

%b20 #&#
\bibitem{ElKaroui-Peng-Quenez}
%
\begin{barticle}[mr]
\bauthor{\bsnm{El~Karoui},~\bfnm{N.}\binits{N.}},
\bauthor{\bsnm{Peng},~\bfnm{S.}\binits{S.}} \AND
\bauthor{\bsnm{Quenez},~\bfnm{M.~C.}\binits{M.~C.}}
(\byear{1997}).
\btitle{Backward stochastic differential equations in finance}.
\bjournal{Math. Finance}
\bvolume{7}
\bpages{1--71}.
\bid{doi={10.1111/1467-9965.00022}, issn={0960-1627}, mr={1434407}}
\bptok{imsref}%
\end{barticle}
%
\endbibitem

%b21 #&#
\bibitem{FlemingSheu}
%
\begin{barticle}[mr]
\bauthor{\bsnm{Fleming},~\bfnm{W.~H.}\binits{W.~H.}} \AND
\bauthor{\bsnm{Sheu},~\bfnm{S.~J.}\binits{S.~J.}}
(\byear{2000}).
\btitle{Risk-sensitive control and an optimal investment model}.
\bjournal{Math. Finance}
\bvolume{10}
\bpages{197--213}.
\bid{doi={10.1111/1467-9965.00089}, issn={0960-1627}, mr={1802598}}
\bptok{imsref}%
\end{barticle}
%
\endbibitem

%b22 #&#
\bibitem{Foellmer-Schied}
%
\begin{bbook}[mr]
\bauthor{\bsnm{F{\"o}llmer},~\bfnm{Hans}\binits{H.}} \AND
\bauthor{\bsnm{Schied},~\bfnm{Alexander}\binits{A.}}
(\byear{2004}).
\btitle{Stochastic Finance: An Introduction in Discrete Time},
\bedition{extended} ed.
\bseries{de Gruyter Studies in Mathematics}
\bvolume{27}.
\bpublisher{de Gruyter}, \baddress{Berlin}.
\bid{doi={10.1515/9783110212075}, mr={2169807}}
\bptok{imsref}%
\end{bbook}
%
\endbibitem

%b23 #&#
\bibitem{Frittelli}
%
\begin{barticle}[mr]
\bauthor{\bsnm{Frittelli},~\bfnm{Marco}\binits{M.}}
(\byear{2000}).
\btitle{The minimal entropy martingale measure and the valuation
problem in
incomplete markets}.
\bjournal{Math. Finance}
\bvolume{10}
\bpages{39--52}.
\bid{doi={10.1111/1467-9965.00079}, issn={0960-1627}, mr={1743972}}
\bptok{imsref}%
\end{barticle}
%
\endbibitem

%%b24 #&#
%%
%(\byear{2002}).
%%

%b25 #&#
\bibitem{Harremoes}
%
\begin{barticle}[mr]
\bauthor{\bsnm{Harremo{\"e}s},~\bfnm{Peter}\binits{P.}}
(\byear{2008}).
\btitle{Some new maximal inequalities}.
\bjournal{Statist. Probab. Lett.}
\bvolume{78}
\bpages{2776--2780}.
\bid{doi={10.1016/j.spl.2008.03.028}, issn={0167-7152}, mr={2465121}}
\bptok{imsref}%
\end{barticle}
%
\endbibitem

%b26 #&#
\bibitem{Hiriart-Urruty}
%
\begin{bbook}[auto:STB|2012/06/01|15:11:24]
\bauthor{\bsnm{Hiriart-Urruty},~\bfnm{J.~B.}\binits{J.~B.}} \AND
\bauthor{\bsnm{Lemar{\'e}chal},~\bfnm{C.}\binits{C.}}
(\byear{2004}).
\btitle{Fundamentals of Convex Analysis}.
\bpublisher{Springer}, \baddress{Berlin}.
\bptok{imsref}%
\end{bbook}
%
\endbibitem

%b27 #&#
\bibitem{Hu-Imkeller-Muller}
%
\begin{barticle}[mr]
\bauthor{\bsnm{Hu},~\bfnm{Ying}\binits{Y.}},
\bauthor{\bsnm{Imkeller},~\bfnm{Peter}\binits{P.}} \AND
\bauthor{\bsnm{M{\"u}ller},~\bfnm{Matthias}\binits{M.}}
(\byear{2005}).
\btitle{Utility maximization in incomplete markets}.
\bjournal{Ann. Appl. Probab.}
\bvolume{15}
\bpages{1691--1712}.
\bid{doi={10.1214/105051605000000188}, issn={1050-5164}, mr={2152241}}
\bptok{imsref}%
\end{barticle}
%
\endbibitem

%b28 #&#
\bibitem{Hu-Schweizer}
%
\begin{bincollection}[mr]
\bauthor{\bsnm{Hu},~\bfnm{Ying}\binits{Y.}} \AND
\bauthor{\bsnm{Schweizer},~\bfnm{Martin}\binits{M.}}
(\byear{2011}).
\btitle{Some new {BSDE} results for an infinite-horizon stochastic control
problem}.
In \bbooktitle{Advanced Mathematical Methods for Finance}
\bpages{367--395}.
\bpublisher{Springer}, \baddress{Heidelberg}.
\bid{doi={10.1007/978-3-642-18412-3_14}, mr={2792087}}
\bptok{imsref}%
\end{bincollection}
%
\endbibitem
%
%%b29 #&#
%%
%(\byear{1994}).
%%

%b30 #&#
\bibitem{Kobylanski00}
%
\begin{barticle}[mr]
\bauthor{\bsnm{Kobylanski},~\bfnm{Magdalena}\binits{M.}}
(\byear{2000}).
\btitle{Backward stochastic differential equations and partial differential
equations with quadratic growth}.
\bjournal{Ann. Probab.}
\bvolume{28}
\bpages{558--602}.
\bid{doi={10.1214/aop/1019160253}, issn={0091-1798}, mr={1782267}}
\bptok{imsref}%
\end{barticle}
%
\endbibitem

%b31 #&#
\bibitem{Lenglart-Lepingle-Pratelli}
%
\begin{bincollection}[mr]
\bauthor{\bsnm{Lenglart},~\bfnm{E.}\binits{E.}},
\bauthor{\bsnm{L{\'e}pingle},~\bfnm{D.}\binits{D.}} \AND
\bauthor{\bsnm{Pratelli},~\bfnm{M.}\binits{M.}}
(\byear{1980}).
\btitle{Pr\'esentation unifi\'ee de certaines in\'egalit\'es de la th\'eorie
des martingales}.
In \bbooktitle{Seminar on {P}robability, {XIV} ({P}aris, 1978/1979)
({F}rench)}.
\bseries{Lecture Notes in Math.}
\bvolume{784}
\bpages{26--52}.
\bpublisher{Springer}, \baddress{Berlin}.
\bid{mr={0580107}}
\bptok{imsref}%
\end{bincollection}
%
\endbibitem

%b32 #&#
\bibitem{LepeltierSanMartin97}
%
\begin{barticle}[mr]
\bauthor{\bsnm{Lepeltier},~\bfnm{J.~P.}\binits{J.~P.}} \AND
\bauthor{\bsnm{San~Martin},~\bfnm{J.}\binits{J.}}
(\byear{1997}).
\btitle{Backward stochastic differential equations with continuous
coefficient}.
\bjournal{Statist. Probab. Lett.}
\bvolume{32}
\bpages{425--430}.
\bid{doi={10.1016/S0167-7152(96)00103-4}, issn={0167-7152}, mr={1602231}}
\bptok{imsref}%
\end{barticle}
%
\endbibitem

%b33 #&#
\bibitem{LepeltierSanMartin98}
%
\begin{barticle}[mr]
\bauthor{\bsnm{Lepeltier},~\bfnm{J.~P.}\binits{J.~P.}} \AND
\bauthor{\bsnm{San~Mart{\'{\i}}n},~\bfnm{J.}\binits{J.}}
(\byear{1998}).
\btitle{Existence for {BSDE} with superlinear-quadratic coefficient}.
\bjournal{Stochastics Stochastics Rep.}
\bvolume{63}
\bpages{227--240}.
\bid{issn={1045-1129}, mr={1658083}}
\bptok{imsref}%
\end{barticle}
%
\endbibitem

%%b34 #&#
%%
%(\byear{2003}).
%measure}.
%%

%b35 #&#
\bibitem{Mania-Schweizer}
%
\begin{barticle}[mr]
\bauthor{\bsnm{Mania},~\bfnm{Michael}\binits{M.}} \AND
\bauthor{\bsnm{Schweizer},~\bfnm{Martin}\binits{M.}}
(\byear{2005}).
\btitle{Dynamic exponential utility indifference valuation}.
\bjournal{Ann. Appl. Probab.}
\bvolume{15}
\bpages{2113--2143}.
\bid{doi={10.1214/105051605000000395}, issn={1050-5164}, mr={2152255}}
\bptok{imsref}%
\end{barticle}
%
\endbibitem

%b36 #&#
\bibitem{ManiaTevzadze}
%
\begin{barticle}[mr]
\bauthor{\bsnm{Mania},~\bfnm{Michael}\binits{M.}} \AND
\bauthor{\bsnm{Tevzadze},~\bfnm{Revaz}\binits{R.}}
(\byear{2006}).
\btitle{An exponential martingale equation}.
\bjournal{Electron. Commun. Probab.}
\bvolume{11}
\bpages{206--216 (electronic)}.
\bid{doi={10.1007/978-3-540-30788-4_24}, issn={1083-589X}, mr={2266711}}
\bptok{imsref}%
\end{barticle}
%
\endbibitem

%%b37 #&#
%%
%(\byear{1976}).
%In \bbooktitle{S\'eminaire de {P}robabilit\'es, {X} ({S}econde Partie:
%{T}h\'eorie des Int\'egrales Stochastiques, {U}niv. {S}trasbourg,
%{S}trasbourg, Ann\'ee Universitaire 1974/1975)}.
%%

%b38 #&#
\bibitem{Mocha-Westray}
%
\begin{bmisc}[auto:STB|2012/06/01|15:11:24]
\bauthor{\bsnm{Mocha},~\bfnm{M.}\binits{M.}} \AND
\bauthor{\bsnm{Westray},~\bfnm{N.}\binits{N.}}
(\byear{2011}).
\bhowpublished{Quadratic semimartingale BSDEs under and exponential moments
condition. Working paper. Available at arXiv:\arxivurl{1101.2582v1}}.
\bptok{imsref}%
\end{bmisc}
%
\endbibitem

%b39 #&#
\bibitem{Morlais1}
%
\begin{barticle}[mr]
\bauthor{\bsnm{Morlais},~\bfnm{Marie-Am{\'e}lie}\binits{M.-A.}}
(\byear{2009}).
\btitle{Quadratic {BSDE}s driven by a continuous martingale and
applications to
the utility maximization problem}.
\bjournal{Finance Stoch.}
\bvolume{13}
\bpages{121--150}.
\bid{doi={10.1007/s00780-008-0079-3}, issn={0949-2984}, mr={2465489}}
\bptok{imsref}%
\end{barticle}
%
\endbibitem

%b40 #&#
\bibitem{Morlais2}
%
\begin{barticle}[auto:STB|2012/06/01|15:11:24]
\bauthor{\bsnm{Morlais},~\bfnm{M.~A.}\binits{M.~A.}}
(\byear{2010}).
\btitle{A new existence result for BSDEs with jumps and
application to
the utility maximization problem}.
\bjournal{Stochastic Process. Appl.}
\bvolume{120}
\bpages{1966--1995}.
\bptok{imsref}%
\end{barticle}
%
\endbibitem

%b41 #&#
\bibitem{Neveu}
%
\begin{bbook}[mr]
\bauthor{\bsnm{Neveu},~\bfnm{Jacques}\binits{J.}}
(\byear{1972}).
\btitle{Martingales \`a Temps Discret}.
\bpublisher{Dunod}, \baddress{Paris}.
\bid{mr={0402914}}
\bptok{imsref}%
\end{bbook}
%
\endbibitem

%b42 #&#
\bibitem{PardouxPeng90}
%
\begin{barticle}[mr]
\bauthor{\bsnm{Pardoux},~\bfnm{{\'E}.}\binits{{\'E}.}} \AND
\bauthor{\bsnm{Peng},~\bfnm{S.~G.}\binits{S.~G.}}
(\byear{1990}).
\btitle{Adapted solution of a backward stochastic differential equation}.
\bjournal{Systems Control Lett.}
\bvolume{14}
\bpages{55--61}.
\bid{doi={10.1016/0167-6911(90)90082-6}, issn={0167-6911}, mr={1037747}}
\bptok{imsref}%
\end{barticle}
%
\endbibitem

%b43 #&#
\bibitem{Protter}
%
\begin{bbook}[mr]
\bauthor{\bsnm{Protter},~\bfnm{Philip~E.}\binits{P.~E.}}
(\byear{2005}).
\btitle{Stochastic Integration and Differential Equations},
\bedition{2nd} ed.
\bseries{Stochastic Modelling and Applied Probability}
\bvolume{21}.
\bpublisher{Springer}, \baddress{Berlin}.
\bid{mr={2273672}}
\bptok{imsref}%
\end{bbook}
%
\endbibitem

%b44 #&#
\bibitem{Rouge-ElKaroui}
%
\begin{barticle}[mr]
\bauthor{\bsnm{Rouge},~\bfnm{Richard}\binits{R.}} \AND
\bauthor{\bsnm{El~Karoui},~\bfnm{Nicole}\binits{N.}}
(\byear{2000}).
\btitle{Pricing via utility maximization and entropy}.
\bjournal{Math. Finance}
\bvolume{10}
\bpages{259--276}.
\bid{doi={10.1111/1467-9965.00093}, issn={0960-1627}, mr={1802922}}
\bptok{imsref}%
\end{barticle}
%
\endbibitem%\

%b45 #&#
\bibitem{Tevzadze}
%
\begin{barticle}[mr]
\bauthor{\bsnm{Tevzadze},~\bfnm{Revaz}\binits{R.}}
(\byear{2008}).
\btitle{Solvability of backward stochastic differential equations with
quadratic growth}.
\bjournal{Stochastic Process. Appl.}
\bvolume{118}
\bpages{503--515}.
\bid{doi={10.1016/j.spa.2007.05.009}, issn={0304-4149}, mr={2389055}}
\bptok{imsref}%
\end{barticle}
%
\endbibitem

\end{thebibliography}
\end{document}